\documentclass[12pt,a4paper]{article}
\usepackage[hebrew,english]{babel}

\usepackage{amsfonts,amssymb}
\usepackage{amsmath}

\numberwithin{equation}{section}

\newtheorem{theorem}{Theorem}

\newtheorem{prop}{Proposition}
\newtheorem{definition}{Definition}
\newtheorem{remark}{Remark}
\newtheorem{corollary}{Corollary}

\newtheorem{lemma}{Lemma}
\newtheorem{example}{Example}

\newcommand{\qed}{\(\Box\)}
\newcommand{\openrm}{\mathrm{(}}
\newcommand{\closerm}{\mathrm{)}}
\newcommand{\dd}{\mathrm{d}}

\def\dim{\hbox{\rm dim}\,}

\def\deg{\hbox{\rm deg}\,}

\numberwithin{lemma}{subsection} \numberwithin{theorem}{subsection}
\numberwithin{prop}{subsection} \numberwithin{remark}{subsection}
\numberwithin{definition}{subsection} \numberwithin{corollary}{subsection}

\begin{document}

\begin{titlepage}

\thispagestyle{empty}

\begin{center}

$ $

\vspace{5cm}\Huge{\bf Ramsey properties of subsets of
\(\mathbb{N}\)}

\vspace{5cm} \large{Thesis submitted for the degree\\
``Doctor of Philosophy"\\}

\vspace{1cm} \large {by\\}

\vspace{1cm} {\Large \bf Alexander Fish}

\vspace{1.5cm}

\vspace{1.5cm} \normalsize{Submitted to the Senate of the Hebrew
University}

\vspace{1cm} \today

\end{center}
\end{titlepage}

\newpage
$ $ \thispagestyle{empty}

\newpage

$ $

\vspace{5cm}{\large This work was carried out under the
supervision of\\}

\hspace{0.5cm}{\large Prof. Hillel Furstenberg}

\thispagestyle{empty}

\newpage
$ $ \thispagestyle{empty}

\newpage

$ $

\vspace{5cm}{\large With my appreciation, I wish to thank my advisor \\}
{\large Prof. Hillel
Furstenberg for his guidance, support \\}
{\large and encouragement during these unforgettable \\}
{\large years of learning and research \\}

\thispagestyle{empty}

\newpage
$ $ \thispagestyle{empty}

\newpage

\tableofcontents
\newpage
$ $ \thispagestyle{empty}
\newpage
\section*{Abstract}
\addcontentsline{toc}{section}{Abstract}
 We  associate
ergodic properties to some subsets of the natural numbers. For any
given family of subsets of the natural numbers one may study the
question of occurrence of certain "algebraic patterns" in every subset in
the family. By "algebraic pattern" we mean a set of solutions of a
system of diophantine equations. In this work we investigate a
concrete family of subsets - WM sets. These sets are characterized by
the property that the dynamical systems associated to such sets are "weakly mixing", and as such they
represent a broad family of randomly constructed  subsets of \(
\mathbb{N} \). We find that certain systems of equations are solvable within every
WM set, and our subject is to learn which systems have this property. We give a complete characterization of linear
diophantine systems which are solvable within every WM set. In
addition we study some non-linear equations and systems of
equations with regard to the question of solvability within every
WM set.
\newpage
$ $ \thispagestyle{empty}
\newpage
\section{Introduction}

The aim of  the dissertation is to  develop Ramsey theory as it
relates to a special family of subsets of the natural numbers,
namely, WM sets.

\subsection{Number theoretic aspects of Ramsey theory and Dynamics}
There are several domains in mathematics where the phenomena of Ramsey theory are encountered. The most classical one
is graph theory.
One of the best examples of Ramsey type theorems in graph theory is Ramsey's theorem:

\noindent
\textit{For every \( k \in \mathbb{N} \) there exists a natural number \( N \) big enough such that
for every coloring into two colors of edges of the complete graph with \( N \) vertices there will exist
 a monochromatic complete
 subgraph with \( k \) vertices.}

\noindent Throughout our work \( \mathbb{N} \) denotes the natural numbers.

\noindent Another theorem of the same spirit, where after a finite coloring of a  structure
we can find a substructure of the same type
at least in one of colors , is van der Waerden theorem:

\noindent
\textit{For every \( r,l \in \mathbb{N}\) there exists \(N(r,l)\in \mathbb{N} \) such that if the integers
\{1,2,\ldots,N(r,l)\} are partitioned into \(r\) sets, one of these contains arithmetic progressions of length
\( l+1\).}

\noindent Note that here \( \{1,2,,\ldots,N(r,l)\}\) may be replaced by any arithmetic
progression of the same length.

\noindent Both Ramsey and van der Waerden theorems may be formulated in the following way:

 \noindent
 \textit{After partitioning into a finite number of
 subsets of a "highly organized" structure (set) we will necessarily find one  subset which contains
 the same substructure.}

 \noindent The difference between the two theorems is in the choice of  "structure".

\noindent The foregoing finite version of van der Waerden theorem is equivalent to the following
 claim
about finite partitions of the natural numbers:

\noindent
\textit{For every partitioning of \(\mathbb{N}\) into a finite number of sets \( C_1, \ldots,C_r \) at least
  one of the subsets contains arbitrarily long arithmetic progressions.}

\noindent
In the thirties of the twentieth century it was conjectured by Erd\"{o}s and Tur\'{a}n  that
the pattern of arbitrarily long arithmetic progressions is not only stable for finite partitions but it necessarily
appears in every subset of the natural numbers with positive upper Banach density. Later
 this conjecture was established by
Szemer\'{e}di, see  \cite{szemeredi}:

\noindent
\textit{The subsets of \( \mathbb{N} \) of positive upper Banach density contain arbitrarily
long arithmetic progressions.}

\noindent The structure of an arithmetic progression of length \( k \) can be viewed as a solution of the following
diophantine system:

\[
  \left\{ \begin{array}{llll} x_2-x_1=x_3-x_2 \\
x_2-x_1 = x_4-x_3 \\
\ldots \\
x_2 - x_1=x_k - x_{k-1}.
\end{array} \right.
\]

\noindent In this work we use extensively  the notion of "algebraic pattern" or, to be more precise,
 we will speak of a subset \( S \) of natural numbers as containing
  some algebraic pattern. The latter means that for some diophantine system
of equations in \( k \)
variables, the set of solutions of the system intersects with \( S^k \).
 Every pattern in this work will be an algebraic pattern.
For example, an arithmetic
progression of length \( k \) is an algebraic pattern.

\noindent There are algebraic patterns which are \textit{regular} for finite
 partitions of \( \mathbb{N} \);
i.e., one of the subsets of the partition necessarily contains the algebraic pattern, but no simple density
condition implies that the pattern will be found.
As an
example of this we present Schur's theorem, \cite{schur}:

\noindent
\textit{For every partitioning of \(\mathbb{N}\) into a finite number of sets \( C_1, \ldots,C_r \) at least
  one of the subsets contains \( x,y,z \) such that x+y=z.}

\noindent It is obvious that  positivity of density for a subset \( S \) is not enough to
ensure existence of a
"Schur pattern" (e.g. \( S =\)odd numbers). In the context of van der Waerden and Schur theorems it will
be appropriate to recall that there is a common generalization of them,
Rado's theorem, which is a complete
characterization of all linear patterns regular for finite partitions. By the word linear
we mean that all
equations in the diophantine system connected to the pattern are linear.

\noindent In the work we are motivated by the following question:

\noindent
\textit{Are there conditions on \( S \subset \mathbb{N} \) more restrictive than positive density
that yield more algebraic patterns?}

\noindent Our way to answer to the question is to add a
condition of "random" behavior (which will be defined rigorously in the next subsection)
to positivity of density of a subset.
A subset which satisfies the foregoing two conditions (is called WM set) will contain
Schur patterns.
Here we would like to give a simple example of "random" behavior.

\noindent We recall that \textit{an infinite \( \{0,1\}\)-valued sequence \( \lambda \) is called a
\textbf{normal sequence} if
every finite binary word \( w \) occurs in \( \lambda \) with a right frequency \(\frac{1}{2^{|w|}}\), where \( |w| \)
is a length of \( w \).} The more familiar notion is that of a normal number \( x \in [0,1]\).
 For every \( x \in [0,1]\), except a countable number of \( x \)'s, there exists
a unique dyadic expansion: \( x = \sum_{i=1}^{\infty} \frac{x_i}{2^i} \, , \, \forall i\,:\,x_i \in \{0,1\}\). Then \( x \) is called a normal number if
and only if  the sequence \( (x_1,x_2,\ldots,x_n,\ldots) \) is a normal sequence.
 To a sequence
\( \lambda=(\lambda_1,\lambda_2,\ldots,\lambda_n,\ldots) \in \{0,1\}^{\mathbb{N}} \) we associate the set
\( B_{\lambda} \subset \mathbb{N} \) by the  rule: \( i \in B_{\lambda} \leftrightarrow \lambda_i=1 \). We
define the notion of a normal set.

\noindent
\textit{A set \( S \subset \mathbb{N}\) is called \textbf{normal} if there exists a normal sequence
\( \lambda \in \{0,1\}^{\mathbb{N}} \) such that \( B_{\lambda} = S \).}

\noindent Normal sets exhibit a non-periodic, "random" behavior.
We remark that every normal set contains Schur patterns.
 We notice that if \( S \) is a normal set then \( S - S \) contains \( \mathbb{N} \).
  Therefore, the equation
 \( z-y=x \) is solvable within every normal set. From the last statement it follows that
 every normal set
 contains Schur patterns.

\noindent We are looking for a possible answer to the foregoing question by using a dynamical approach.
All aforementioned theorems have dynamical equivalent formulations. For our question the most relevant
theorem is Szemer\'{e}di's  theorem. Furstenberg has shown that Szemer\'{e}di's  theorem is equivalent
to the phenomenon of multiple recurrence valid for general volume preserving
dynamical systems which can be established  by purely dynamical techniques
 (see \cite{furst4}).

\noindent   In this context Furstenberg formulates a correspondence principle for subsets
of the natural numbers of positive upper Banach density:

\noindent
\textit{Given a set \( E \subset \mathbb{N} \) with \( d^*(E) > 0 \) (\( E \) of positive upper Banach density) there exists
a probability measure preserving system \( (X,\mathbb{B},\mu,T)\) and a set \( A \in \mathbb{B}, \mu(A) = d^*(E)\),
 such that
for any \( k \in \mathbb{N}\) and any \( n_1,\ldots,n_k \in \mathbb{Z} \) one has:
\[
  d^*(E\cap(E-n_1)\cap \ldots \cap(E-n_k)) \geq \mu(A\cap T^{-n_1}A \cap \ldots \cap T^{-n_k}A).
\]}

\noindent By this correspondence principle in order to prove Szemer\'{e}di's theorem it is sufficient
to establish the following multiple recurrence theorem which is proved purely dynamically in \cite{furst4}.

\noindent
\textit{For any probability measure preserving system \( (X,\mathbb{B},\mu,T)\), a set \( A \in \mathbb{B}, \mu(A) > 0\)
and any \( k \in \mathbb{N} \) there exists \( n \in \mathbb{N} \) such that
\( \mu(A\cap T^{-n}A \cap T^{-2n}A \cap \ldots \cap T^{-(k-1)n}A ) > 0\).
}

 \noindent  The basic idea of the correspondence principle is that a set of positive density can be viewed
 more or less (there are some technicalities) as return times of generic points of  ergodic systems
  to a set of
 positive measure. If a dynamical system will be even more "random" (for example weakly mixing or mixing) then
 we expect to find that within a set of return times to a set of positive measure
 one can find a greater variety of algebraic patterns.

 \noindent Our approach is to deal with the sets of integers that are the return times of a generic point of
 weakly mixing system to a set of positive measure. Such subsets of \( \mathbb{N} \) we call WM sets.
 We formalize this in the next section.

\subsection{Generic points and WM sets}\label{intro} To define formally the
main object of this work we need the notions of
measure preserving systems and of  generic points.

\begin{definition}
Let \( X \) be a compact metric space, \( \mathbb{B} \) be the Borel
\( \sigma\)-algebra on \( X \), let \( T:X \rightarrow X \) be a
measurable map and \( \mu \) a probability measure on \(
\mathbb{B} \). A quadruple \( (X,\mathbb{B},\mu,T) \) is called a
\textbf{measure preserving system} if for every \( B \in
\mathbb{B} \) we have \( \mu(T^{-1}B) = \mu(B) \).
\end{definition}

\noindent For a compact metric space \( X \) we denote by \( C(X)
\) the space of continuous functions on \( X \) with the uniform
norm.
\begin{definition} Let \( (X,\mathbb {B},\mu ,T) \) be a measure
preserving system. A point \( \xi \in X \) is called
\textbf{generic} if for any \( f\in C(X) \) we have
\begin{equation}
\label{limit} \lim _{N\rightarrow \infty }\frac{1}{N}\sum
^{N-1}_{n=0}f(T^{n}\xi )=\int _{X}f(x)d\mu (x).
\end{equation}
\end{definition}


\noindent We can now give an alternative definition of a normal
set which is purely dynamical. A set \( S \) is normal if and only
if the sequence \( 1_S \in \{0,1\}^{\mathbb{N}} \) is a generic
point of the measure preserving system \((
\{0,1\}^{\mathbb{N}},\mathbb{B},T,\mu)\), where \( \mathbb{B} \)
is Borel \( \sigma \)-algebra on the topological space \(
\{0,1\}^{\mathbb{N}} \) which is endowed with the Tychonoff topology,
\( T \) is the shift to the left, \( \mu \) is the product measure of
\(\mu_i\)'s where  \( \mu_i(0)=\mu_i(1) = \frac{1}{2} \).
Thus, the system \(( \{0,1\}^{\mathbb{N}},\mathbb{B},T,\mu)\)
is the Bernoulli \( (\frac{1}{2}, \frac{1}{2})\) system and,
in particular, it is a mixing dynamical system.

\noindent   The notion of a WM set generalizes that of a normal set, where the role
played by a Bernoulli dynamical system is taken over by dynamical systems  of more general
 character.
\newline
Let \( \xi (n) \) be any \( \{0,1\}- \)valued sequence. There is a
natural dynamical system \( (X_{\xi },T) \) connected to the
sequence \( \xi \):

\noindent On the foregoing compact space \( \Omega
=\{0,1\}^{\mathbb {N}} \) which is endowed with the Tychonoff topology,
we define a continuous map \( T:\Omega \longrightarrow \Omega  \)
by shifting all the elements of a sequence to left, namely, \(
(T\omega )_{n}=\omega _{n+1} \). Now for any \( \xi \) in \(
\Omega \) we define \( X_{\xi } \) to be \(
\overline{(T^{n}\xi )_{n\in \mathbb {N}}} \subset \Omega\). 

\noindent Let \( A \) be a subset of \( \mathbb {N} \). Choose \(
\xi =1_{A} \) and assume that for an appropriate measure \( \mu \),
the point \( \xi \) is generic for \( (X_{\xi },\mathbb {B},\mu
,T) \). Now we attach to the set \( A \) dynamical properties
associated with the system \( (X_{\xi },\mathbb {B},\mu ,T) \).
\newline
For example, \( A \) is called weakly mixing (respectively - totally
ergodic) if the measure preserving system \( (X_{\xi },\mathbb
{B},\mu ,T) \) is weakly mixing (respectively - totally ergodic).
\newline
We recall the latter two notions of ergodic theory.
\begin{definition}
A measure preserving system \( (X,\mathbb{B},\mu,T) \) is called
\textbf{ergodic} if every \( A \in \mathbb{B} \) which is
invariant under \(T\), i.e. \( T^{-1}(A) = A \), satisfies \(
\mu(A) = 0 \) or \(1 \).
\newline
A measure preserving system \( (X,\mathbb{B},\mu,T) \)  is called
\textbf{totally ergodic} if for every \( n \in \mathbb{N} \) the
system \( (X,\mathbb{B},\mu,T^n) \) is ergodic.
\newline
A measure preserving system \( (X,\mathbb{B},\mu,T) \) is called
\textbf{weakly mixing} if the system \( (X\times X,\mathbb{B}_{X
\times X},\mu \times \mu,T \times T) \) is ergodic.
\end{definition}
\noindent Let \(\mathcal{P}\) denote  some dynamical property of a measure preserving system.
We can attach the property \( \mathcal{P} \) to a subset of the natural numbers by the following:
\begin{definition}
A subset \( S \subset \mathbb{N} \) is \( \mathcal{P} \) \( \Leftrightarrow \) \( 1_S \) is generic for
measure preserving system \( (X_{1_S},\mathbb{B},\mu,T)\) which has property \( \mathcal{P}\).
\end{definition}

\noindent Finally,
we would like to deal with subsets of \( \mathbb{N} \) which
may have a rich structure, i.e. may be expected to exhibit many algebraic
patterns. Therefore, we restrict ourselves to the case of weakly
mixing subsets of \( \mathbb{N}\) of positive density (the density
of every weakly mixing set exists!). For completeness we define
the density of a subset of \( \mathbb{N} \).
\begin{definition}
Let \( S \subset \mathbb{N} \). If the limit of \( \frac{1}{N} \sum_{n=1}^{N} 1_S(n) \) exists
 as \( N \rightarrow \infty \) we call it the \textbf{density} of \( S \) and denote by \( d(S) \).
\end{definition}

\begin{definition}
A subset \( S \subset \mathbb{N} \) is called a \textbf{WM set} if
\( S \) is weakly mixing and the density of \( S \) is positive.
That is to say, \( 1_S \) is a generic point of the weakly mixing
system \( (X_{1_S},\mathbb{B},\mu ,T) \) and \( d(S) > 0 \).
\end{definition}

\noindent We could equally well speak of strongly mixing sets, but
for our purposes, weak mixing will be adequate.

\subsection{Examples of combinatorial properties of WM sets}
\label{combin_properties}
\noindent We would like to list  basic combinatorial/Ramsey properties of WM
sets. To do this we recall the definitions of two basic notions in
ergodic Ramsey theory and combinatorial number theory.

\begin{definition}
\label{Poincare_def}
A set \( S \subset \mathbb{N} \) is called a \textbf{Poincar\'{e}
set} if for every measure preserving  system \( (X,\Sigma,T,\mu)\)
(not necessarily topological system) and every \( A \in \Sigma \)
with \( \mu(A) > 0 \) there exists \( n \in  S \) such that \(
\mu(A \cap T^{-n}A) > 0 \).
\end{definition}
This can be reformulated in purely combinatorial terms.
First, we recall the notion of upper Banach density for a subset of the natural numbers.
\begin{definition}
Let \( E \subset \mathbb{N} \). Upper Banach density of \( E \), \( d^*(E) \) is the following quantity
\[
  d^*(E) = \limsup_{b_n-a_n \rightarrow \infty} \frac{|E \cap \{a_n,\ldots,b_n\}|}{b_n-a_n+1}.
\]
\end{definition}

\noindent \textnormal{By Furstenberg's correspondence principle, for a set \( S \) to be Poincar\'{e} is equivalent to the following:}
\newline
\textit{For every subset \( E \) of positive upper Banach density there exists \( s \in S \) with \( d^*(E\cap(E-s))>0 \).}
\newline
In fact, a milder condition is sufficient: \textit{a set \( S \) is Poincar\'{e} if and only if
for every \( E \) of positive
upper Banach density there exists \( s \in S \) such that
\(
E \cap (E-s) \neq \emptyset.
\);} that is to say that \( s \) is a difference of two numbers in \( E \).

\noindent The last property is called \( 1 \)-recurrence. Sometimes in the literature a
 Poincar\'{e}
set is called \( 1 \)-recurrent set.

\noindent The next notion is taken from combinatorial number
theory and may be viewed as a generalization of an infinite
arithmetic progression around \( 0 \).

\begin{definition}
 A set  \( S \subset \mathbb{N} \) is called \textbf{IP-set} if there
exists an infinite sequence of  natural numbers \(
\{p_1,p_2,\ldots,p_n,\ldots \}\) (not necessarily different) such
that \[ S = \{ p_{i_1}+\ldots+p_{i_k} \, | \, i_1 < i_2 \ldots <
i_k \, , \, k \in \mathbb{N}\}.
\]
\end{definition}
We recall the definition of IP*-set.
\begin{definition}
A set \( S \) is called an IP*-set if for every IP-set \( E \) we have \( E \cap S \neq \emptyset \).
\end{definition}

\noindent The next two results which will be proved in Chapter
\ref{sect_basic_Ramsey_properties_WM_sets} give a first evidence
of the richness of algebraic patterns which occur in every WM set.

\begin{theorem}
\label{Poincare}
Every WM set is a \( 1 \)-recurrent set. (\( \Rightarrow \) Poincar\'{e} set)
\end{theorem}

\begin{theorem} \label{IP}
Every WM set contains an IP-set.
\end{theorem}
\begin{corollary} A non-trivial WM set (which has density less than \( 1 \)) is never an IP*-set.
\end{corollary}

\noindent One of the reasons to choose WM sets as an object of our
research  and not sets which satisfy weaker conditions, for
example, totally ergodic sets, is the fact that the foregoing theorems
don't hold for totally ergodic sets. The following is an example
of a totally ergodic set which is neither a Poincar\'{e} set nor
contains an IP-set.

\begin{example}
Let \( \alpha \not \in \mathbb{Q} \) and denote by \( S \) the
following subset of \(\mathbb{N}\)
\[
S = \left\{n \in \mathbb{N} \, | \, \alpha n \, mod \, 1 \in
\left[\frac{2}{5},\frac{3}{5}\right] \right\}.
\]
Then \( S \) is a totally ergodic set of positive density which is
not Poincar\'{e} set and for every \( x,y \in S \) we have \( x+y
\not \in S \).
\end{example}
\begin{proof}
We start from the last statement which is easily proven. Namely,
if \( \xi,\eta \in [\frac{2}{5},\frac{3}{5}] \) then \( (\xi + \eta) \mod{1}  \not \in
[\frac{2}{5},\frac{3}{5}] \). It follows that
if \(
x,y \in S \) then \( (x+y) \alpha \mod 1 \not \in
\left[\frac{2}{5},\frac{3}{5}\right] \) and therefore \( x+y \not
\in S \). This implies \( S \) contains no IP-set.
\newline
\( S \) is not a Poincar\'{e} set, as we see by checking the recurrence condition
 of definition \ref{Poincare_def} for the system \(
(\mathbb{T},\mathbb{B},S_{\alpha},\lambda) \), where \(
\mathbb{T}\) is the one dimensional torus, \( S_{\alpha}(x) = x +
\alpha \), \( \lambda \) is lebesgue measure and the subset \(A =
\left[0,\frac{1}{5} \right]\) is of  measure \( \frac{1}{5} \).
Then obviously for every \( s \in S \) we have \( \lambda(A
\cap S_{\alpha}^{-s}A) = 0 \).
\newline
To show that \( S \) is a totally ergodic set we note that \( S \)
consists of  return times to the set \( I =
\left[\frac{2}{5},\frac{3}{5}\right] \) of the point zero within
the aforementioned measure preserving system \(
(\mathbb{T},\mathbb{B},S_{\alpha},\lambda) \). Consider
the space of \(\{0,1\}\)-sequences \( X_{1_S}\), and consider the
characteristic function \( \chi \in C(X_{1_S}) \) of a cylinder
\[ C_{i_1,\ldots,i_k}^{j_1,\ldots,j_k} = \{ \omega \in \{0,1\}^{\infty} \, |
\,\omega_{i_l} = j_l \, , \, \forall 1 \leq l \leq k\}. \]

\noindent We have
\[
  \frac{1}{N} \sum_{n=1}^{N} \chi(T^n 1_S) = \frac{1}{N} \sum_{n=1}^{N} \phi_{j_1}(n + i_1)\ldots \phi_{j_k}(n+i_k)=
\]
\[
  \frac{1}{N} \sum_{n=1}^{N} T^n f(0) \rightarrow_{N \rightarrow \infty} \int_{\mathbb{T}} f(x) d \lambda(x),
\]
where \( \phi_1(n) = 1_I(\alpha n)\),\( \phi_0(n) = 1 - 1_I(\alpha
n)\) and \( f(x) = \prod_{l=1}^k T^{i_l}\phi_{j_l}(x) \). Since
the linear space of characteristic functions on cylinders is dense
in \( C(X_{1_S}) \) we conclude that the point \( 1_S \) is  a
generic point in \( X_{1_S} \) for a measure which is obtained as
a projection of lebesgue measure in
\((\mathbb{T},\mathbb{B},S_{\alpha},\lambda) \).Thus our system
is a factor of a totally ergodic system
\((\mathbb{T},\mathbb{B},S_{\alpha},\lambda) \); therefore it is
itself totally ergodic.

\hspace{12cm} \qed
\end{proof}

\noindent The concept of  WM sets is new. It relies on properties of the
corresponding point \( 1_S \) within a measure preserving dynamical system. There is a
 concept of "good" subsets in the context topological dynamics (they contain many algebraic
patterns) which is defined by H. Furstenberg (see \cite{furst3});
namely, \textit{central sets}. In order to define these sets we define a uniformly recurrent point
in a topological dynamical system.
\begin{definition}
Let \( (X,T) \) be a topological dynamical system, i.e., \( X \)
is a metric compact space and \( T:X \rightarrow X \) is a
continuous transformation. A point \( x_0 \in X \) is called
\textbf{uniformly recurrent} if for any open set \( U \), such
that \( x_0 \in U \), the set \( \{ n \in \mathbb{N} | T^n x_0 \in
U \} \) is syndetic (a set with bounded gaps).
\end{definition}

\begin{definition}
Let \( (X,T) \) be a topological dynamical system. Denote by \( d
\) a metric on \( X \). Two points \( x,y \in X \) will be called
\textbf{proximal} if there exists an increasing sequence \(
\{n_k\}\) such that \( \lim_{k \rightarrow \infty}
d(T^{n_k}x,T^{n_k} y) = 0\).
\end{definition}
\begin{definition}
A set \( S \subset \mathbb{N} \) is called \textbf{central} if
there exists a topological dynamical system  \( (X,T) \),
 a uniformly recurrent point \( x_0 \in X \), a point \( x \in X \) which is proximal to \( x_0 \) and
 a neighborhood \( U \) of \( x_0 \) such that \( S = \{ n \in \mathbb{N} | T^n x \in U \} \).
\end{definition}

\noindent In the section \ref{central_sets} we prove the
incomparability of central and WM sets.
\begin{theorem}
\label{central}
There exists a WM set which does not contain  a central set.
\end{theorem}
\begin{remark}
\textnormal{The opposite direction is easy; for example, we
could take the set of even numbers.}
\end{remark}

\subsection{Main results}\label{mainresults}
\subsubsection{Solvability of linear diophantine systems within WM
sets}\label{sub_sect_linear_dioph_syst}

 \noindent We have succeeded to give a complete
characterization of those linear systems of diophantine equations
which are solvable within every WM set.
\begin{theorem}
\label{main_thm_lin}
 Let \( B \in \mathbb{Q}^{t \times k} \) and \(
\vec{d} \in \mathbb{Q}^t \). The system of linear equations
\begin{equation}
\label{lin_eq}
 B\vec{x}=\vec{d}
\end{equation}
 is solvable within  every WM set \( \Leftrightarrow \)
there exist three vectors \( \vec{x_1} = (a_1,a_2,\ldots,a_k)^t ,\,
\vec{x_2} = (b_1,b_2,\ldots,b_k)^t, \, \vec{f} = \{f_1,f_2,\ldots,f_k\}^t \in \mathbb{N}^k \),
 disjoint
sets \( E,F_1,\ldots,F_l \subset \{1,2,\ldots,k\}\), \( E \cup F_1 \cup \ldots F_l =
\{1,2,\ldots,k\}\), such that:
\newline
a\(\closerm\) for every \( i,j \in E,\, i \neq j \)
\[
\det \left( \begin{array}{cc}
                    a_i & b_i \\
                    a_j & b_j  \\
            \end{array}
    \right) \neq 0.
\]
\newline
b\(\closerm\) for every \( p \in \{1,\ldots,l\}\) there exist \( c_1^p,c_2^p \in \mathbb{N} \), such that
for every \( i \in F_p \) we have \( a_i = c_1^p \, , \, b_i = c_2^p \)
and for every \(j \in \{1,\ldots,k\} \setminus F_p \) we require
\[
\det \left( \begin{array}{cc}
                    a_j & b_j \\
                    c_1^p & c_2^p  \\
            \end{array}
    \right) \neq 0.
\]
\newline
c\(\closerm\) The vector \( \vec{f} \) is constant on all indices from the same \( F_p \) with \( p \in \{1,2,\ldots,l\}\), namely,
\[ \forall p \in \{1,2,\ldots,l\} \, \forall i \in F_p \, : \, f_i = f^p, \]
where \( f^p \in \mathbb{Z} \).
d\(\closerm\) The affine space of solutions of the system \(
B\vec{x}=\vec{d} \) contains \[
\{ n\vec{x_1} + m\vec{x_2} + \vec{f} \, | \, n,m \in \mathbb{N} \}.
\]

%
\end{theorem}

\noindent A proof of this theorem is in  section \ref{lin_section} of the thesis.

\noindent As we will show in proposition \ref{Rado_analog}, it will follow from theorem \ref{main_thm_lin}
 that every linear algebraic pattern which is
regular for finite partitions (Rado theorem gives a complete characterization of such patterns)
occurs in every WM set. In the context of Rado patterns and WM sets it is
important to recall that by Furstenberg's theorem (see \cite{furst3}) every central set contains all
Rado patterns. By theorem \ref{central} it follows that we can't prove that every WM set contains all
Rado patterns by use of Furstenberg's theorem.

\subsubsection{An additive analog of polynomial multiple
recurrence for WM Sets }\label{additive_polynomial}

\noindent  A natural generalization of the theorem of
Szemer\'{e}di is the seminal theorem of Bergelson and Leibman
about polynomial multiple recurrence \cite{berg_leibman}. If we
rephrase this theorem combinatorially it states that for every \(
k \) polynomials which are essentially distinct (i.e., no two
differ by a constant) \( p_1,\ldots,p_k \) with positive leading
coefficients and \( p_1(0) = p_2(0)=\ldots=p_k(0) = 0 \), and for
every subset \( A \) of the natural numbers of  positive upper
Banach density,  there exists \( n \in \mathbb{N} \)
 such that \( \{x,x+p_1(n),\ldots,x+p_k(n)\} \in A^{k+1} \). The
latter means that the system of equations
\[
\left\{ \begin{array}{lll} y_1-x=p_1(n) \\
 \ldots \\
y_k - x = p_k(n)
\end{array}
\right.
\]
is solvable in every set of positive upper Banach density for some \( n \in \mathbb{N} \). For \( A \) a WM set
we can use Bergelson's PET theorem (see \cite{berg_pet}) and to obtain the same result without the restriction that
all polynomials have zero free coefficient. If
additionally we require that \( n \in A \) then  we can
use the IP-polynomial Szemer\'{e}di theorem of Bergelson, Furstenberg
and McCutcheon (see \cite{berg_furst_mc}) and the fact that every
WM set contains an IP-set to establish that the previous system is
solvable within every WM set \( A \) and  \( n \in
A \) provided \( p_i(0) = 0 \, , \, \forall i: \, 1 \leq i \leq k \).
 It is very natural question to try to establish the
analogous result for the "additive" system which is obtained from
the last one by replacing all minuses by pluses.
\begin{equation}
\label{additive_system} \left\{ \begin{array}{llll} x+y_1=p_1(z) \\
x+y_2 = p_2(z) \\
\ldots \\
x+y_k=p_k(z)
\end{array} \right.
\end{equation}
Of course, in the case of the additive system we can't expect that
there exists a solution within every set of positive upper density
(there are a lot of examples of periodic sets that contain no
solution for the equation \( x + y = n^2 \); i.e., the set \( 5 \mathbb{N} + 1 \)). On the other hand, we
would expect that for some such systems there exists a solution
within every WM set, where congruence conditions do not form an obstruction. We can obtain the following 
characterization of solvability of system (\ref{additive_system})
within every WM set.
\begin{theorem}
\label{main_thm}
 For every \( k \in \mathbb{N} \) the system
(\ref{additive_system}) is solvable within every WM set if 
\( \deg (p_1) = \deg(p_2) = \ldots = \deg (p_k) \), the difference of every
two polynomials is a non-constant polynomial and all leading coefficients of \( p_1,\ldots,p_k\)
are positive.
\end{theorem}

\noindent There is an easy case which shows the necessity of some
restrictions on the degrees of the polynomials;  namely, when in
the system (\ref{additive_system}) there are two polynomials with
degrees which differ by at least two.
\begin{remark}
\label{remark_dif_at_least_two}
 \textnormal{If in the system
(\ref{additive_system}) there are two polynomials with degrees
which differ by at least two, then there exist WM sets within which
the system (\ref{additive_system}) is unsolvable.}
\end{remark}
\begin{proof}
We take an arbitrary WM set \( A \); then removing a set of
density zero from \( A \) leads again to a WM set. In particular, we can exclude from \( A \) all
solutions of the system
(\ref{additive_system}) by removing a set of density zero.
Namely, if \( \deg{p_1} \leq \deg{p_2} - 2 \) then replace \(
A \) by
\[
A' = A \setminus \left( \bigcup_{n \in \mathbb{N}}
[p_2(n)-p_1(n),p_2(n)] \right)
\]
which is again a WM set. (For sufficiently large \( n \) the polynomials \(
p_1(n), p_2(n) \) are monotone.) Within \( A' \) the system
(\ref{additive_system}) is unsolvable.

\hspace{12cm} \qed
\end{proof}

\subsubsection{The equation $xy = z$ and normal sets}
\label{normal_sets_xyz}

\noindent We recall the notion of a normal set.
\newline
We have the natural bijection between infinite binary \(
\{0,1\}\)-sequences and subsets of \( \mathbb{N} \), namely for
any sequence \( \lambda \) we associate the subset \( B_{\lambda}
= \{i| \lambda_i = 1 \} \).
\begin{definition}
A set \( B \subset \mathbb{N} \) is called \textbf{normal} if the
infinite binary sequence \( \lambda \) which corresponds to \( B
\) (i.e. \(B_{\lambda} = B\)) is normal.
\end{definition}
In the section \ref{normal_section} we prove the following result.
\begin{theorem}
There exist normal sets within which the equation \( xy=z \) is
unsolvable.
\end{theorem}

\noindent Our proof is non-constructive and we do not know an explicit example.

\noindent
On the other hand the equation \( xy=z^2\) is solvable in any
normal set, and in fact:
\begin{theorem}
\label{th1_pos} Let \( A \subset \mathbb{N} \) be a WM set. Then
there exist \( x,y,z \in A \) (\( x \not = y \)) such that \( xy =
z^2 \).
\end{theorem}

\noindent For normal sets we can also show the following
\begin{theorem}
\label{th2_pos}
 Let \( A \subset \mathbb{N} \) be an arbitrary
normal set. Then there exist \( x,y,u,v \in A \) such that \( x^2 +
y^2 = square \) and \( u^2 - v^2 = square \).
\end{theorem}

\noindent This result holds for WM sets as well provided their
density exceeds \( \frac{1}{3}\).

\subsection{Structure of the thesis}\label{sub_sect_struct_of_thesis}

\noindent The thesis consists of \(5\) sections and an Appendix. The first section is an introduction to the subject of the thesis,
namely WM sets, and a formulation of main results. In the second section we prove basic combinatorial properties of WM sets,
which rely on 1-recurrence of WM sets. In addition we show that the notions of central sets and of WM sets are incomparable.
In the third section we give a proof of the theorem which  characterizes all linear diophantine systems which are solvable
within every WM set. In the fourth section we prove that the system (\ref{additive_system}) is solvable within every WM set
 if 
 all the polynomials are essentially distinct, have the same degree and have positive leading coefficients. The section \(5\)
  is devoted to non-linear
 equations.
 In particular, we show the existence of a normal set for which the (non-linear) equation
 \( xy=z \) has no solutions with \( x,y,z\) in the set.
In Appendix we collected some technical lemmas which are used in more than one section.
In particular we formulate and prove the van der Corput lemma which will be used
on several occasions.

\newpage
\section{Basic combinatorial properties of WM sets}\label{sect_basic_Ramsey_properties_WM_sets}

\subsection{Every WM set is a Poincar\'{e} set}\label{subsect_WM_is_Poincare}

\noindent In this section we prove the following theorem.

\noindent
\textbf{Theorem \ref{Poincare}} \textit{Every WM set is a \( 1 \)-recurrent set. (\( \Rightarrow \) Poincar\'{e} set)}

\noindent To prove theorem \ref{Poincare} we note that it is sufficient by the ergodic decomposition
theorem to show recurrence of a WM set for every ergodic system. We show the following

\begin{prop}
Let \( S \) be a WM set. Then for every ergodic measure preserving system
\( (X,\Sigma,\mu,T) \) and every \( A \in \Sigma \) with \( \mu(A) > 0 \) there exists
\( s \in S \) such that \( \mu(A \cap T^{-s} A) > 0 \).
\end{prop}
\begin{proof}
We make  use of spectral theory. Namely, by spectral theory for
the unitary operator  \( U: L^2(X,\Sigma,\mu) \rightarrow L^2(X,\Sigma,\mu) \) which is
defined by \( U f = f \circ T \) and the function \( 1_A \in
L^2(X,\mu) \) there exists a spectral measure \( \omega_{1_A} \) (
we denote it simply \( \omega \)) on \( \mathbb{T} \) (the
spectrum of \( U \)) such that for every \( n \in \mathbb{N} \)
we have
\[
  <1_A,T^n 1_A > = \int_{[0,1]} e^{2\pi i x n} d\omega(x).
\]
Let \( \mathbb{T} \) denotes \(1\)-dimensional torus and
for every \( \alpha: \, 0 \leq \alpha < 1 \) let
\newline
\( S_{\alpha}(x) \doteq x + \alpha (\mod{1})\).
For every \( \alpha \in (0,1) \) consider the Kronecker system
\( (K, S_{\alpha}) \), where \( K = \overline{\{S_{\alpha}^n(0)\}_{n=0}^{\infty}}\).
For \( \alpha \not \in \mathbb{Q} \)  this system is
a factor of the system \((\mathbb{T},\mathbb{B},S_{\alpha},\lambda)\), defined in
\( \S \) \ref{combin_properties}, and in any case
the Kronecker system is
disjoint from  the weak-mixing
system \( (\overline{\{T^n 1_S\}_{n=1}^{\infty}}, \mathbb{B}, \mu,
T ) \). We use the theorem of Furstenberg:
\newline
\textit{If the measure preserving systems \( (X,\mathbb{B}_X,\mu,T_X)\) and \( (Y,\mathbb{B}_Y,\nu,T_Y)\) are disjoint,
\( x \in X \) and \( y \in Y \) are generic then \( (x,y) \in X \times Y \) is generic for the system
\((X \times Y, \mathbb{B}_X \times \mathbb{B}_Y, \mu \times \nu, T_X \times T_Y) \) (see \cite{furst2}).}

\noindent Applying this to the pair \((0,1_S) \in \mathbb{T} \times \overline{\{T^n 1_S\}}\)
 we obtain
\begin{equation}
\label{Poincare_spectral_eq}
  \lim_{N \rightarrow \infty} \frac{1}{N} \sum_{n=1}^N 1_S(n) e^{2 \pi i \alpha n} = 0.
\end{equation}
Therefore, by use of Lebesgue dominated convergence theorem from (\ref{Poincare_spectral_eq}) we have
\[
  \frac{1}{N} \sum_{n=1}^N 1_S(n) <1_A,T^n 1_A > \rightarrow_{N \rightarrow \infty} d(S) \omega(0).
\]
But \( \omega(0) \) is represented in terms of integral over \( 1_A \) by the following
\[
  \omega(0) = \int_{[0,1]} \lim_{N \rightarrow \infty} \frac{1}{N} \sum_{n=1}^{N} e^{2\pi i x n} d\omega(x) =
\]
\[
  \lim_{N \rightarrow \infty} \frac{1}{N} \sum_{n=1}^{N} \int_{[0,1]} e^{2\pi i x n} d\omega(x) =
  \lim_{N \rightarrow \infty} \frac{1}{N} \sum_{n=1}^{N} <1_A,T^n 1_A > =
  <1_A, \lim_{N \rightarrow \infty} \frac{1}{N} \sum_{n=1}^{N} T^n 1_A > =
\]
\[
  \left( \int_X 1_A(x) d \mu(x) \right)^2 = \mu(A)^2.
\]
We have used ergodicity of the system \( X \) in the last step.

\noindent Finally, we get
\[
  \frac{1}{N} \sum_{n=1}^N 1_S(n) <1_A,T^n 1_A > \rightarrow_{N \rightarrow \infty} d(S) \mu(A)^2 > 0.
\]
Since the inner product \( <1_A,T^n 1_A>= \mu(A \cap T^{-n} A)\), we conclude that there exists \( s \in S \), such that \( \mu(A \cap T^{-s} A) > 0 \).

\hspace{12cm} \qed
\end{proof}

\subsection{Every WM set contains an IP set}\label{subsect_WM_contains_IP}

\noindent To prove that every WM set contains an IP set we use
theorem \ref{Poincare}.

\noindent
\begin{proof} (of  theorem \ref{IP})
\newline
Let \( S \) be an arbitrary WM set. By use of theorem
\ref{Poincare} we conclude that there exists \( s \in S \) such
that \( (S - s) \cap S \) has positive density. We shall see that it is again a WM set.
\newline
To prove the last statement we define in the weak-mixing measure
preserving space \( (X = \overline{\{T^n 1_S\}_{n=1}^{\infty}},
\mathbb{B},\mu,T) \) the set \( A = \{ x \in X \, | \, (x)_0 =
1 \} \). Then \( \mu(A) = d(S) > 0 \) (by use of genericity of \(
1_S \) in \( X \)) and by using recurrence of the set \( S \) we
obtain that there exists \( s \in S \), such that \( \mu(A \cap
T^{-s}A)
> 0 \). By genericity of the point \( 1_S \in X \) it follows
that \( \mu(A \cap T^{-s}A) = d(S \cap (S-s)) \).
The map \( \phi: X \rightarrow \{0,1\}^{\mathbb{N}} \) defined by \(\phi(x) = y \), where
\( y(n) = x(n)x(n+s) \) takes \( X \) to a closed shift invariant set \( Y \) in
\(\{0,1\}^{\mathbb{N}}\) and \(\phi(1_S)=1_{S \cap S-s}\).
\newline
Therefore the point \( 1_{S \cap (S-s)} \in \{0,1\}^{\mathbb{N}} \) is a
generic point of \( (Y,T) \) which is again a weak-mixing measure preserving system.
Here we get weak-mixing  because the resulting system is
a factor of the system \( (X,\mathbb{B},\mu,T)\).

\noindent The next stage of our proof is to define inductively  an IP set in \( S \).
\newline
Let \( s_1 \in S \), such that \( S \cap (S - s_1) \) is again a WM set.
\newline
If we denote by \( S_1 = S \cap (S-s_1) \) (a WM set) then we define \( s_2 \in S_1 \), such that
\( S_1 \cap (S_1 - s_2) \) is again a WM set. Note that if \( s_3 \in S_1 \cap (S_1-s_2) \)
then \( s_3+s_1,s_3+s_2,s_3+s_1+s_2 \in S \).
\newline If we have defined \( s_1,\ldots,s_n \) and a WM set \( S_n \)
we define the element \( s_{n+1} \in S_n \) and a WM set \( S_{n+1} \) by the following: there
exists an element \( s_{n+1} \in S_n \), such that \( S_n \cap (S_n - s_{n+1}) \) is again a WM set, which we denote by
\(S_{n+1}\).

\noindent In this way we have defined an infinite sequence \( \{s_1,s_2,\ldots,s_n,\ldots\} \subset S\). It is a consequence
of the construction of the sequence that every finite sum of its elements is again in \( S \).
\newline
Therefore we have found an IP set within an arbitrary WM set.

\hspace{12cm} \qed

\end{proof}

\subsection{Incomparability of Central and WM sets}\label{central_sets}

\noindent We will use a  variant of Rohlin's lemma in ergodic theory.
\begin{lemma} \( \label{lem1} \) Let \( (X, \mathbb{B}, \mu, T) \) be an ergodic non-periodic invertible
measure preserving system (m.p.s.). Then for any \( n > 1 \) there
exists \( C \in \mathbb{B} \) with \( \frac{1}{n} \leq \mu(C) \leq
\frac{1}{n-1} \),
 such that
\( X = \bigcup_{i=0}^{n} T^i C \), where the equality is up to a
measure zero set and every point \( x \in \bigcup_{i=0}^{n} T^i
C \) returns to \( C \) by at most \( n+1 \) iterations of \(  T
\).
\end{lemma}
\begin{proof}
Let us fix \( n > 1 \). For any \( B \in \mathbb{B} \) with \( \mu(B)
> 0 \) let us build up Kakutani's tower, by the following
procedure.
\newline
Denote by \( B_k = \{ x \in B | r_B(x) = k \} \), where
\( r_B(x) = \min_{ i \geq 1} \{ i | T^i x \in B \} \). By the Poincar\'{e} recurrence theorem
 we have
 \( B = \bigcup_{ i=1 }^{\infty} B_i \).
Then the following family of sets
 \( B_1 , ( B_2 \cup T B_2 ) , \ldots , (B_k \cup \ldots \cup T^{k-1} B_k)  , \ldots \)
is called Kakutani's tower with base \( B \).
\newline
Obviously by ergodicity it follows \( X = \bigcup_{k=1}^{\infty} \bigcup_{i=0}^{k-1} T^i B_k \)
(where the equality is up to a set of measure zero) and the union is measurably disjoint
(this means that the intersection of any two sets
from the union is of measure zero).
\newline
We introduce \( C = \bigcup_{m=0}^{\infty} \bigcup_{ k = m(n+1) +
1}^{\infty} T^{m(n+1)} B_k \). Then \(  X = \bigcup_{i=0}^{n} T^i
C \), from which we get the estimation \( \frac{1}{n} \leq \mu(C)
\), and any point in \( C \) returns back to \( C \) by at most \(
n+1 \) iterations of \( T\). On the other hand if we denote by \(
B' \) the higher layer of \( C \), then obviously we have
\[
 \mu(C \setminus B') \leq \frac{1}{n}.
\]
In addition we have \( \mu(B') = \mu(B) \). Therefore we
get
\[
\mu(C) \leq \frac{1}{n} + \mu(B).
\]
\noindent Now let us choose \( B \in \mathbb{B} \) such that \( 0
< \mu(B) < \frac{1}{n-1} - \frac{1}{n} = \frac{1}{n(n-1)} \) (it
can be done because our m.p.s. is non-periodic and is non-atomic).
Finally with appropriate choice of \( B \) we get the desired
result.

\hspace{12cm} \qed
\end{proof}

\begin{remark} \( \label{rem1} \)
We don't have to assume that the system is non-periodic; it is sufficient that for any \( \varepsilon > 0 \)
there exists \( B \in \mathbb{B} \) with \( 0 < \mu(B) \leq \varepsilon \).
\end{remark}

\begin{remark} \( \label{rem2} \)
Every weak-mixing m.p.s. is non-periodic.
\end{remark}
\begin{definition}
A sequence \( \omega \in \{ 0 , 1 \}^{\mathbb{N}} \) is said to
satisfy \textbf{property \( (l,L) \)} if for any \( l > 0 \) there
exists \( L
> 0 \) such that for every \( k \geq 0 \) the block \( (\omega_k,
\ldots, \omega_{k+L-1} ) \) contains at least one subblock of \( l
\) successive zeros.
\end{definition}
We attach the notion of centrality to \( \{0,1\}\)-valued sequences as well.
\begin{definition}
A sequence \( \lambda \in \{0,1\}^{\mathbb{N}} \) is called \textbf{central} if the set
\( B_{\lambda} = \{ i \in \mathbb{N} \, | \, \lambda_i = 1 \}\) is central.
\end{definition}
\begin{lemma} \( \label{lem2} \)
Let \( \omega \in \{ 0 , 1 \}^{\mathbb{N}} \) be a sequence satisfying
property \( (l,L) \),
then \( \omega \) is not a central sequence.
\end{lemma}
\begin{proof}
First of all we prove that if \( \omega' \) is proximal to  \(
\omega \) and \( \omega' \)  is uniformly recurrent then \(
\omega' \)  is the zeros sequence. For, let us assume that \(
\omega' \) has aforementioned properties and take a block \( (
\omega'_0, \omega'_1, \ldots , \omega'_{m-1} ) \). Then there
exists \( l \) such that this block is contained in any continuous
subblock of length \( l \) of \( \omega' \). But  \( \omega \) is
proximal to \( \omega' \), thus for any \( L > 0 \) there exists
\( n \geq 0 \) such that \( \{\omega\}_{n}^{n+L-1} =
\{\omega'\}_{n}^{n+L-1} \). Let us choose \( L = L(l) \) such that
any subblock of \( \omega \) of length \( L \) contains \( l \)
successive zeros.
\newline
As a result of our choices, the block \( ( \omega'_0, \omega'_1, \ldots , \omega'_{m-1} ) \) is a subblock of
any subblock of \( \{\omega\}_{n}^{n+L-1} \) of length \( l \), in particular the zeros block of
length \( l \).
\newline
Thus we have proved that \( \omega \) is proximal to only one uniform recurrent sequence, namely the zeros
sequence.
\newline
Suppose, contrary to the assertion of the lemma that \( \omega \) is a central sequence. Then
there exists \( V \), a neighborhood of the zeros sequence, such
that \(  \omega_n = 1\) iff \( T^n \omega \in V \). \( V \)
contains an open set which contains the zeros sequence, thus there
exists \( l \geq 1 \), such that any word \( x \) that begins with
\( l \) zeros is inside \( V \) (\( l \) can not be zero, because
then \( \omega \) is ones's sequence which does not satisfy \(
(l,L) \) property); \( \omega \) is proximal to the zeros
sequence, thus there exists \( k \geq 0 \) such that  \(
\{\omega\}_{k}^{k+l-1} \) is the zeros block. But in this case \(
T^k \omega \in V \) and therefore \( \omega_k = 1 \).
\newline
Thus we have got a contradiction to our assumption that \( \omega \) is central.

\hspace{12cm} \qed
\end{proof}

\begin{theorem} \( \label{thm} \)
Let \( (X, \mathbb{B}, \mu, T) \) be a weak mixing invertible m.p.s.. Then there exist a symbolic weak mixing system
\( (Y, \mathbb{B}, \nu, T) \) and \( y_0 \in \{ 0 , 1 \}^{\mathbb{N}} \) a generic point of \( Y \) with
\( 0 < d(y_0) < 1 \) (where \( d \) is the density of ones) such that
every sequence \( y \leq y_0 \)
(for every \( n \) we have
\(y(n) \leq y_0(n)) \) is not a central sequence  and \( Y \) is a factor of \( X \).
\end{theorem}
\begin{proof}
By remark \ref{rem2}  the system \( (X, \mathbb{B}, \mu, T) \)
satisfies all the requirements of lemma \ref{lem1}  and therefore for
any \( n > 1 \) there exists \( C_n \in \mathbb{B} \) such that \(
X = \bigcup_{i=0}^{n} T^i C_n \) with \( \frac{1}{n} \leq \mu
(C_n) \leq \frac{1}{n-1} \),
 and every
point inside \( C_n \) returns back to \( C_n \) by at most \( n+1 \) iterations of \( T \).
\newline
Now we construct \( A \in \mathbb{B} \) of a positive measure which is bounded by a predefined
number \( \alpha \) by the
following procedure.
\newline
Let us choose \( \{ L_l \} \) a sequence of positive natural
numbers (for every \( l\) we assume that \( L_l \geq 2 \)) such
that \( \sum_{l=1}^{\infty} \frac{l}{L_l-1} = \alpha \). Then for
any \( L_l \) let us take \( C_{L_l} \in \mathbb{B} \) as above
with \( X = \bigcup_{i=0}^{L_l} T^i C_{L_l} \) and , finally, take
\[
A = \bigcup_{l=1}^{\infty}( \bigcup_{ i=0 }^{ l-1 } T^i C_{L_l})
\]
The following estimations on the measure of \( A \) are obvious:
\[
 \frac{1}{L_1} \leq \mu(A) \leq \sum_{l=1}^{\infty} l \mu(C_{L_l}) \leq \sum_{l=1}^{\infty} \frac{l}{L_l-1} = \alpha
\]
Let \( \phi_A : X \rightarrow \{ 0 ,1 \}^{\mathbb{N}} \) to be
defined by \( \omega = \phi_A(x) \) iff \( \omega(n) = 1 - 1_A(T^n
x) \). Then obviously \( \phi_A \) is measurable and \( \phi_A
\circ T = T \circ \phi_A \) (where \( T \) on the right hand is
the usual shift transformation).
 Let us define \( Y \risingdotseq \phi_A (X) \), then
\( ( Y , \mathbb{B}_Y , ( \phi_A )_* \mu, T ) \) is a m.p.s. and a factor of \( X \), thus is a weak mixing system
(for any \( B \in \mathbb{B}_Y \) we define \( (\phi_A )_* \mu \risingdotseq \mu( \phi_A^{-1} (B) )\).
\newline
 Let us denote by \( X' \) the following subset of \( X \)
\[
X' \risingdotseq  \bigcap_{l=1}^{\infty}( \bigcup_{ i=0 }^{ L_l + 1 } T^i C_{L_l})
\]
It is obvious that \( \mu (X') = \mu (X) \). Let \( G \subset X \)
be the set of generic points in \( X \) (\( X \) is a compact
metric space and therefore the notion of a generic point is well
defined). By the ergodic theorem \( \mu (G) = \mu(X) \), and
therefore \( \mu (G \cap X') = \mu(X) \) and thus the measure of
\( \phi_A(G \cap X' ) \) in Y is equal to the measure of whole \(
Y \).
 But \( (Y,T) \) is ergodic (even weak mixing)
therefore almost every point of \( \phi_A(G \cap X') \)  is
generic. Choose \( y_0 \in \phi_A(G \cap X') \) to be generic in
\( Y \). Then there exists \( x_0 \in G \cap X' \) such that \(
y_0 = \phi_A (x_0) \) (by using the ergodic theorem once again, we
can add one more condition on \( x_0 \), namely \( \frac{1}{N}
\sum_{n=1}^{N} 1_A(T^n x_0) \rightarrow \mu(A) \) [it should be
done because \( 1_A \) might be a non continuous function]). It is
obvious that \( d(y_0) = \mu ( A^c ) \) and thus \( d(y_0) \geq 1
- \alpha \).
\newline
By the choice of \( A \) and \( x_0 \) it follows that \( y_0 \) satisfies  \( (l,L) \)  property.
Every sequence \( y \leq y_0 \) is again satisfies \((l,L)\) property and, thus,
by lemma \ref{lem2}, it  follows
that  \( y \) is not a central sequence.

\hspace{12cm} \qed
\end{proof}

\begin{remark} \( \label{rem2_1} \)
For any \( 0 < \alpha  < 1 \) we can construct \( y_0 \) in
the formulation of the theorem with \( d(y_0) \geq \alpha \).
\end{remark}


\noindent By combining theorem \ref{thm} with remark \ref{rem2_1} we obtain the following statement.
\begin{theorem} \( \label{th_wm_not_c} \)
For any \( 0 < \alpha < 1 \) there exists \( A \in \mathbb{N} \)
 a WM set with \(  d(y_0) \geq \alpha \)
which is not central and such that no subset  \( B \) of \( A \) is a
central sequence.
\end{theorem}

%

\newpage
 \section{Solvability of linear diophantine equations
within WM sets}\label{lin_section}

\subsection{Proof of Sufficiency}\label{lin_proof_suff}

\noindent We restate the main result of this section which was formulated in section \ref{sub_sect_linear_dioph_syst}.
\newline
\textbf{Theorem \ref{main_thm_lin}} \textit{
 Let \( B \in \mathbb{Q}^{r \times k} \) and \(
\vec{d} \in \mathbb{Q}^r \). The system of linear equations
\begin{equation}
\label{lin_eq}
 B\vec{x}=\vec{d}
\end{equation}
 is solvable within  every WM set \( \Leftrightarrow \)
there exist two vectors \( \vec{x_1} = (a_1,a_2,\ldots,a_k)^t ,\,
\vec{x_2} = (b_1,b_2,\ldots,b_k)^t \in \mathbb{N}^k \), disjoint
sets \( E,F_1,\ldots,F_l \subset \{1,2,\ldots,k\}\), \( E \cup F_1 \cup \ldots F_l =
\{1,2,\ldots,k\}\), such that:
\newline
a\(\closerm\) for every \( i,j \in E,\, i \neq j \)
\[
\det \left( \begin{array}{cc}
                    a_i & b_i \\
                    a_j & b_j  \\
            \end{array}
    \right) \neq 0.
\]
\newline
b\(\closerm\) for every \( p \in \{1,\ldots,l\}\) there exist \( c_1^p,c_2^p \in \mathbb{N} \), such that
for every \( i \in F_p \) we have \( a_i = c_1^p \, , \, b_i = c_2^p \)
and for every \(j \in \{1,\ldots,k\} \setminus F_p \) we require
\[
\det \left( \begin{array}{cc}
                    a_j & b_j \\
                    c_1^p & c_2^p  \\
            \end{array}
    \right) \neq 0.
\]
\newline
c\(\closerm\) There exist \( f^1,\ldots,f^l \in \mathbb{Z} \)  such that setting \( f_i = f^p \) for
 \(  p \in \{1,\ldots,l\} \) and \( i \in F_p \), then 
  the affine space of solutions of a system \(
B\vec{x}=\vec{d} \) contains \[
 \{ (a_1 n + b_1 m + f_1, \ldots , a_k n + b_k m + f_k)^t
\hspace{0.1 in} | \hspace{0.1 in} n,m \in \mathbb{N} \}. \]
}

\noindent
\textbf{Notation:} \textit{We introduce the scalar product of two vectors \( v,w \) of the
length \( N \) as follows:
\[
<v,w>_N \doteq \frac{1}{N} \sum_{n=1}^N v(n)w(n).
\]
We denote by \( L^2(N) \) the Hilbert space of all real vectors of
the length \( N \) with the aforementioned scalar product.
\newline
  We define: \(
\parallel{w}\parallel_N^2 \doteq <w,w>_N \).}

\noindent First we state the following proposition which is a very useful tool
in the proof of the sufficiency of the conditions of theorem \ref{main_thm_lin}.
\begin{prop}
\label{main_prop}
 Let \( A_i \subset \mathbb{N} \) \( \rm( \) \( 1 \leq i \leq k \)\( \rm) \) be  WM sets. Let \newline
 \(\xi_i(n) \doteq 1_{A_i}(n) - \dd(A_i) \), where \( \dd(A_i) \) denotes
density of \( A_i \). Suppose there are  \( (a_1,b_1),(a_2,b_2),
\ldots, (a_k,b_k) \in ( \mathbb{Z} \setminus \{0\} )^2 \), such
that \( a_i
> 0, \, 1 \leq i \leq k \), and for every \( i \neq j \)
\begin{displaymath}
\det \left( \begin{array}{cc}
a_i & b_i \\
a_j & b_j  \\
\end{array}
\right) \neq 0.
\end{displaymath}
 Then for every \( \varepsilon > 0 \) there exists \( M(\varepsilon)
 \in \mathbb{N}
 \), such that for every \( M \geq M(\varepsilon) \) there exists
 \( N(M,\varepsilon) \in \mathbb{N} \), such that for every \( N \geq
 N(M,\varepsilon)\)
 \begin{displaymath}
 \left\Vert{w}\right\Vert_N < \varepsilon,
 \end{displaymath}
 where \( w(n) \doteq \frac{1}{M} \sum_{m=1}^M \xi_1(a_1n + b_1 m) \xi_2(a_2 n + b_2 m)
 \ldots \xi_k(a_k n + b_k m) \) for every  \( n=1,2,\ldots,N\).
\end{prop}

\noindent Since the proof of proposition \ref{main_prop} involves many
technical details, first we show how our main result follows from
it. Afterwards we state and prove all the lemmas necessary for a proof of
proposition \ref{main_prop} and define all the required concepts.
\newline
We will need an easy consequence of proposition \ref{main_prop}.
\begin{corollary}
\label{corol_main_prop} Let \( A \) be a WM set. Let \( k \in
\mathbb{N} \), suppose \( (a_1,b_1),(a_2,b_2), \ldots, (a_k,b_k) \in
( \mathbb{Z} \setminus \{0\} )^2 \) satisfy all requirements of
proposition \ref{main_prop} and suppose \( f_1,\ldots,f_k \in
\mathbb{Z} \). Then for every \( \delta > 0 \) there exists \(
M(\delta)\) such that \( \forall \, M \geq M(\delta) \) there exists
\( N(M,\delta) \) such that \( \forall \, N \geq N(M,\delta)\) we
have
\[
\left| \Vert v \Vert_N - d^k(A) \right| < \delta,
\]
where \( v(n) \doteq \frac{1}{M} \sum_{m=1}^M 1_A(a_1n + b_1 m +f_1)
1_A(a_2 n + b_2 m +f_2)
 \ldots 1_A(a_k n + b_k m + f_k) \) for every  \( n=1,2,\ldots,N\).
\end{corollary}
\begin{proof}
We can write \( v(n) \) in the following form:
\[
v(n) =  \frac{1}{M} \sum_{m=1}^M (\xi(a_1n + b_1 m ) +d(A))
(\xi(a_2 n + b_2 m ) +d(A)) \ldots  (\xi(a_k n + b_k m )+d(A)),
\]
for every \( n=1,2,\ldots,N\). We again introduce normalized WM
sequences \( \xi_i(n) = \xi(n+f_i) \). Then by use of triangular
inequality and proposition \ref{main_prop} it follows that for big
enough \( M \) and \(N \) (which depends on \(M \)) \( \Vert v
\Vert_N \) is as close as we wish to \( d^k(A) \). The latter
finishes the proof.

 \hspace{12cm} \qed
\end{proof}

\noindent
\begin{proof} (of the theorem \ref{main_thm_lin}, \( \Lleftarrow \))
\newline
By corollary \ref{corol_main_prop} it follows that the vector \( v \) defined by 
\[
v(n) \doteq \frac{1}{M} \sum_{m=1}^M 1_A(a_1n + b_1 m +f_1)
1_A(a_2 n + b_2 m +f_2)
 \ldots 1_A(a_k n + b_k m + f_k),
\]
for every  \( n=1,2,\ldots,N\) is not identically zero for big enough \(M \) and \( N \). The latter is possible only if for some \(n,m \in \mathbb{N} \) we have
\[
(a_1n + b_1 m +f_1, a_2 n + b_2 m +f_2, \ldots, a_k n + b_k m + f_k) \in A^k.
\]

\hspace{12cm} \qed
\end{proof}

\noindent Now we state and prove all the claims that are required in
order to prove proposition \ref{main_prop}.



\begin{definition}
Let \( \xi \) be a WM-sequence of zero average. The autocorrelation
function of\, \( \xi \) of the length \( j \in \mathbb{N} \) with
the shifts
\(\{\{r_1,i_1\},\{r_2,i_2\},\ldots,\{r_j,i_j\}\}\)\newline
 $\openrm$all
shifts are integers$\closerm$ is the sequence \(
\psi^j_{\{r_1,i_1\},\{r_2,i_2\},\ldots,\{r_j,i_j\}}\) which is
defined as follows: for \( j > 1 \)
\[
\psi^j_{\{r_1,i_1\},\{r_2,i_2\},\ldots,\{r_j,i_j\}}(n)=
\]
\[
\psi^{j-1}_{\{r_1,i_1\},\{r_2,i_2\},\ldots,\{r_{j-1},i_{j-1}\}}(n+r_j)
\psi^{j-1}_{\{r_1,i_1\},\{r_2,i_2\},\ldots,\{r_{j-1},i_{j-1}\}}(n+
r_j+ i_j),
\]
for \( j = 1 \) the autocorrelation function is defined as
\[
\psi^1_{\{r_1,i_1\}}(n) = \xi(n+r_1)\xi(n+r_1 + i_1).
\]
\end{definition}
\begin{remark}
For any sequence \( \psi \) we define \(\psi(-n) = 0 \) for every \(
n \in \mathbb{N} \).
\end{remark}
\begin{lemma}
\label{autocorrelation1} Let \( \xi \) be a WM-sequence of zero
average and suppose \( \varepsilon,\delta > 0 \). Then for every
\( j \geq 1 \), \(\{c_1,c_2,\ldots,c_j\} \in (\mathbb{Z} \setminus
\{0\})^j \) and  \( \{r_1,r_2,\ldots,r_j\} \in (\mathbb{Z})^j \)
there exists \( I = I(\varepsilon,\delta,c_1,\ldots,c_n) \), such
that there exists a set \( S \subset [-I,I]^j \) of density at
least \( 1- \delta \) and there exists \( N(I,\varepsilon) \in
\mathbb{N} \), such that for every \( N \geq N(I,\varepsilon)\)
there exists \( L(N,I,\varepsilon) \) such that for every \( L
\geq L(N,I,\varepsilon) \)
\[
  \frac{1}{L} \sum_{l=1}^L  \left( \frac{1}{N}  \sum_{n=1}^N
\psi^j_{\{r_1,c_1 i_1\},\{r_2,c_2 i_2\},\ldots,\{r_j,c_j i_j\}}(l+ b
n)
 \right)^2 < \varepsilon,
\]
for every \( \{i_1,i_2,\ldots,i_j\} \in S \).
\end{lemma}
\begin{proof}
We note that it is sufficient to prove the lemma in the case \( c_1
= c_2 = \ldots = c_j = 1 \), since if the average of nonnegative
numbers over a complete lattice is small, then the average over a
sublattice of a fixed positive density is also small.
\newline
Recall that \( \xi \in X_{\xi} \doteq \overline{\{T^n \xi\}_{n=0}^{\infty}} \subset  supp({\xi})^{\mathbb{N}} \),
where \( T \) is a usual shift to the left on the dynamical system \( supp({\xi})^{\mathbb{N}} \), and
by the assumption that \( \xi \) is a WM-sequence of zero average it follows that \( \xi \) is a generic
point of the weak-mixing system \( (X_{\xi},\mathbb{B_{X_{\xi}}},\mu,T) \) and the function
\( f \, : \, f(\omega) \doteq \omega_0 \) has zero integral.
\newline
We define functions \(
g_{\{r_1,i_1\},\{r_2,i_2\},\ldots,\{r_j,i_j\}} \) on \( X_{\xi} \)
inductively. Let \( g_{\emptyset} \doteq f \).
Define\newline \(
g_{\{r_1,i_1\},\ldots,\{r_{j-1},i_{j-1}\},\{r_j,i_j\}} \doteq
T^{r_j}\left( g_{\{r_1,i_1\},\ldots,\{r_{j-1},i_{j-1}\}}T^{i_j}
g_{\{r_1,i_1\},\ldots,\{r_{j-1},i_{j-1}\}} \right)\).
\newline
\noindent Define the functions \(
g^*_{\{r_1,i_1\},\{r_2,i_2\},\ldots,\{r_j,i_j\}} = \prod_{\epsilon
\in V_j^{*}} f \circ T^{r_1 + \ldots + r_j + \epsilon_1 i_1 + \ldots
\epsilon_k i_j}\), where \( V_j \) is a \( j\)-dimensional discrete
cube \( \{0,1\}^j\) and \( V_j^{*} \) is the whole \(
j\)-dimensional discrete cube except the zero point. (Note that \( g
= ( T^{r_1 + \ldots + r_j} \circ f )g^{*} \), where we have omitted
subscripts.)
\newline
The following has been proven by Host and Kra in \cite{host-kra}
(theorem 13.1):
\newline \textit{ Let \( (X,\mu,T)\) be an ergodic system. Given an
integer \( k \) and \( 2^k \) bounded functions \( f_{\epsilon} \)
on \( X \), \(\epsilon \in V_k\) , the
functions \[
 \prod_{i=1}^k \frac{1}{N_i - M_i} \sum_{n \in [M_1,N_1)
\times \ldots [M_k,N_k)} \prod_{\epsilon \in V_k^{*}}
f_{\epsilon}\circ T^{\epsilon_1 n_1 + \ldots \epsilon_k n_k}
\]
converge in \( L^2(\mu)\) to the limit function \[
 \mathbb{E}\left(
\bigotimes_{\epsilon \in V_k^{*}} f_{\epsilon} | \tau^{[k]^*}
\right) (x),
\]
when \( N_1 - M_1,\ldots,N_k-M_k\) tend to \( +\infty\). The \(
\sigma\)-algebra \( \tau^{[k]^*} \) is identified with the so-called characteristic factor \(
Z_{k-1}(X) \).}
\newline
The characteristic factors \( Z_{k}(X) \) are defined for arbitrary ergodic systems, and what is important
for our purposes is that  in our case of the weak-mixing system \( X_{\xi} \), all the factors \(
Z_{k-1}(X_{\xi})\) are trivial.

\noindent  Therefore, the limit function
in our case will be a constant. By integrating the limit function
we obtain that this constant is equal to \( \prod_{\epsilon \in
V_k^*} \int_{X_{\xi}} f_{\epsilon} d\mu \).

\noindent From the theorem of Host and Kra, applied to the
weak-mixing system \( X_{\xi} \times X_{\xi} \) and the functions \(
f_{\epsilon}(x) = f \otimes f \) for every \( \epsilon \in V_k \),
 we obtain for every  Folner
sequence \( \{F_n\} \) in \( \mathbb{N}^j \) that an average over the multi-index \(\{i_1,\ldots,i_j\}\)
 of \(
g^*_{\{r_1,i_1\},\{r_2,i_2\},\ldots,\{r_j,i_j\}} \otimes
g^*_{\{r_1,i_1\},\{r_2,i_2\},\ldots,\{r_j,i_j\}}\) on \( F_n \)'s
converges to zero (the integral of \( f \otimes f \) is zero). If we
would take another Folner sequence \(  \{G_n\}\) in \( \mathbb{N}^j
\) then for the same \( \{r_1,\ldots,r_j\}\) the closeness of an
average of \( g^*_{\{r_1,i_1\},\{r_2,i_2\},\ldots,\{r_j,i_j\}}
\otimes g^*_{\{r_1,i_1\},\{r_2,i_2\},\ldots,\{r_j,i_j\}}\) on \( G_n
\) to zero depends only the size of the box \( G_n \). Namely, if
all edges of a box are big enough then the aforementioned average is
small.
\newline
\noindent As a result we have
\newline
 \textit{For every \( \varepsilon
> 0 \), \( j \in \mathbb{N} \) and every fixed \(
\{r_1,r_2,\ldots,r_j\} \in \mathbb{N}^j\), there exists a subset \(
R \subset \mathbb{N}^j \) with lower density equal to one, such
that
\begin{equation}
\label{eq:one_lemma311}
 \left( \int_{X_{\xi}}
g_{\{r_1,i_1\},\ldots,\{r_{j-1},i_{j-1}\},\{r_j,i_j\}} d\mu
\right)^2 < \varepsilon,
\end{equation}
 for every \( \{i_1,i_2,\ldots,i_j\} \in R \).}

\noindent We note that \(
\psi^j_{\{r_1,i_1\},\{r_2,i_2\},\ldots,\{r_j,i_j\}}(l+ b n) =
g_{\{r_1,i_1\},\{r_2,i_2\},\ldots,\{r_j,i_j\}}\left( T^{l+bn} \xi
\right) \).
\newline
The definition of the sequences \( \psi^j \) implies
\[
\lim_{L \rightarrow \infty} \frac{1}{L} \sum_{l=1}^L \left(
\frac{1}{N} \sum_{n=1}^N
\psi^j_{\{r_1,i_1\},\{r_2,i_2\},\ldots,\{r_j,i_j\}}(l+bn) \right)^2
\]
\[
= \lim_{L \rightarrow \infty} \frac{1}{L} \sum_{l=1}^L \left(
\frac{1}{N} \sum_{n=1}^N \psi^j_{\{\pm r_1,\pm i_1\},\{\pm r_2,\pm
i_2\},\ldots,\{\pm r_j,\pm i_j\}}(l \pm bn) \right)^2.
\]
Therefore, in order to prove the Lemma \ref{autocorrelation1}\, it
is sufficient to show the following:
\newline
\textit{For every \( \varepsilon,\delta > 0 \) and for a priori
chosen \( r_1,r_2,\ldots,r_j,b \in \mathbb{N} \)
 there exists  \(
I(\varepsilon,\delta) \in \mathbb{N} \), such that for every \( I
\geq I(\varepsilon,\delta) \) there exists a subset \( S \subset
[1,I]^j \) of density at least \( 1 - \delta \) \( \openrm
\)namely, we have \( \frac{|S \cap [1,I]^j|}{I^j} \geq 1-\delta
\)\( \closerm \) and there exists \( N(I,\varepsilon) \in
\mathbb{N} \), such that for every \( N \geq N(I,\varepsilon)\)
there exists \( L(N,I,\varepsilon) \in \mathbb{N} \) such that for
every \( L \geq L(N,I,\varepsilon) \) the following holds for
every \( \{i_1,i_2,\ldots,i_j\} \in S \):
\[
  \frac{1}{L} \sum_{l=1}^L  \left( \frac{1}{N}  \sum_{n=1}^N
\psi^j_{\{r_1,i_1\},\{r_2,i_2\},\ldots,\{r_j,i_j\}}(l+ b n)
 \right)^2 < \varepsilon.
\]
}
\newline
Assume that \( r_1,r_2,\ldots,r_j,b \in \mathbb{N} \). Continuity of
the function \newline \(
g_{\{r_1,i_1\},\{r_2,i_2\},\ldots,\{r_j,i_j\}} \) and the genericity
of the point \( \xi \in X_{\xi} \) yields
\[
\lim_{L \rightarrow \infty} \frac{1}{L} \sum_{l=1}^L \left(
\frac{1}{N} \sum_{n=1}^N
\psi^j_{\{r_1,i_1\},\{r_2,i_2\},\ldots,\{r_j,i_j\}}(l+bn) \right)^2
\]
\[
= \lim_{L \rightarrow \infty} \frac{1}{L} \sum_{l=1}^L \left(
\frac{1}{N} \sum_{n=1}^N T^{bn}
g_{\{r_1,i_1\},\{r_2,i_2\},\ldots,\{r_j,i_j\}}\left( T^{l} \xi
\right)\right)^2
\]
\begin{equation}
\label{eq_lemm_2}
 = \int_{X_{\xi}} \left( \frac{1}{N} \sum_{n=1}^N T^{bn}
g_{\{r_1,i_1\},\{r_2,i_2\},\ldots,\{r_j,i_j\}} \right)^2 d \mu.
\end{equation}
By combining the ergodic theorem, applied to the weak-mixing system
\( (X_{\xi},\mathbb{B},\mu,T^b) \), with disjointness of any
weak-mixing system from the cyclic system on \( b \) elements we note
that
\begin{equation}
\label{imp_eq_22}
 \frac{1}{N} \sum_{n=1}^N T^{bn}
g_{\{r_1,i_1\},\{r_2,i_2\},\ldots,\{r_j,i_j\}} \rightarrow_{N
\rightarrow \infty}^{L^2(X_{\xi})} \int_{X_{\xi}}
g_{\{r_1,i_1\},\{r_2,i_2\},\ldots,\{r_j,i_j\}} d \mu.
\end{equation}
From (\ref{eq:one_lemma311}) there exists \( I(\varepsilon,\delta) \in \mathbb{N} \) big enough,
such that for every \( I \geq I(\varepsilon,\delta) \) there
exists a set \( S \subset [1,I]^j \) of density at least \(
1-\delta \) such that
\[
 \left( \int_{X_{\xi}}
g_{\{r_1,i_1\},\{r_2,i_2\},\ldots,\{r_j,i_j\}} d \mu \right)^2 <
\frac{\varepsilon}{4},
\]
for all  \( \{i_1,i_2,\ldots,i_j\} \in S \).
\newline
From equation (\ref{imp_eq_22}) follows that there exists \(
N(I,\varepsilon) \in \mathbb{N} \), such that for every \( N \geq
N(I,\varepsilon) \) we have
\[
\int_{X_{\xi}}\left( \frac{1}{N} \sum_{n=1}^N T^{bn}
g_{\{r_1,i_1\},\{r_2,i_2\},\ldots,\{r_j,i_j\}} \right)^2 <
\frac{\varepsilon}{2},
\]
for all \( \{i_1,i_2,\ldots,i_j\} \in S \).
\newline
Finally, equation (\ref{eq_lemm_2}) implies that there exists \(
L(N,I,\varepsilon) \in \mathbb{N} \), such that for every \( L \geq
L(N,I,\varepsilon) \) we obtain
\[
 \frac{1}{L} \sum_{l=1}^L \left(
\frac{1}{N} \sum_{n=1}^N
\psi^j_{\{r_1,i_1\},\{r_2,i_2\},\ldots,\{r_j,i_j\}}(l+bn) \right)^2
< \varepsilon,
\]
for all \( \{i_1,i_2,\ldots,i_j\} \in S \).

\hspace{12cm} \qed
\end{proof}

\noindent The following lemma is a generalization of the previous
lemma for a product of several autocorrelation functions.
\begin{lemma}
\label{autocorrelation2}
 Let \(
\psi^{1,j}_{\{r_1^1,i_1\},\{r_2^1,i_2\},\ldots,\{r_j^1,i_j\}},
\ldots ,
\psi^{k,j}_{\{r_1^k,i_1\},\{r_2^k,i_2\},\ldots,\{r_j^k,i_j\}} \) be
autocorrelation functions of length \( j \) of WM-sequences \(
\xi_1, \ldots, \xi_k \) of zero average, \(
\{c_1^1,\ldots,c_j^1,\ldots, c_1^k,\ldots,c_j^k\} \in (\mathbb{Z}
\setminus \{0\})^{jk}\) and \( \varepsilon,\delta
> 0 \). Suppose\newline \( (a_1,b_1),(a_2,b_2), \ldots, (a_k,b_k) \in
\mathbb{Z}^2 \), such that \( a_i > 0, b_i \neq 0 \),\( 1\leq
i\leq{k}\) and for every \( i \neq j \)
\begin{displaymath}
\det \left( \begin{array}{cc}
a_i & b_i \\
a_j & b_j  \\
\end{array}
\right) \neq 0.
\end{displaymath}
Then there exists \(I(\varepsilon,\delta) \in \mathbb{N} \), such
that for every \( I \geq I(\varepsilon,\delta) \) there exist \(
S \subset [-I,I]^j \) of density at least \( 1- \delta \), \(
M(I,\varepsilon) \in \mathbb{N} \), such that for every \( M \geq
M(I,\varepsilon) \) there exists \( X(M,I,\varepsilon) \in
\mathbb{N}\), such that for every \( X \geq X(M,I,\varepsilon) \)
\[
 \frac{1}{X} \sum_{x=1}^X ( \frac{1}{M} \sum_{m=1}^M
\psi^{1,j}_{\{r_1^1,c_1^1 i_1\},\{r_2^1,c_2^1
i_2\},\ldots,\{r_j^1,c_j^1 i_j\}}(a_1x+ b_1 m) \ldots\]
\[
\psi^{k,j}_{\{r_1^k,c_1^k i_1\},\{r_2^k,c_2^k
i_2\},\ldots,\{r_j^k,c_j^k i_j\}}(a_kx+ b_k m))^2 < \varepsilon,
\]
for every \( \{i_1,i_2,\ldots,i_j\} \in S \).
\end{lemma}

\noindent \begin{proof}
The proof is by induction on \( k \). The case \( k=1 \) (and
arbitrary \( j \)) follows from the Lemma \ref{autocorrelation1} and
the Proposition \ref{wm_prop}.
\newline
Suppose that the statement holds for \( k-1 \).
\newline
Denote by \[
 v_m(x) \doteq \psi^{1,j}_{\{r_1^1,c_1^1
i_1\},\ldots,\{r_j^1,c_j^1 i_j\}}(a_1x+ b_1 m) \ldots
\psi^{k,j}_{\{r_1^k,c_1^k i_1\},\ldots,\{r_j^k,c_j^k i_j\}}(a_kx+
b_k m). \] The van der Corput lemma (lemma \ref{vdrCorput} of the
appendix)
 implies that it is
sufficient to show the existence of  \(
\mathbb{I}(\varepsilon,\delta) \in \mathbb{N} \), such that for
every \( \mathbb{I} \geq \mathbb{I}(\varepsilon,\delta) \) there
exists a set \( S \subset [-\mathbb{I},\mathbb{I}]^j \) of density
at least \( 1-\delta \) and there exists \(
I(\varepsilon,\mathbb{I}) \) big enough (\(
I(\varepsilon,\mathbb{I}) \geq I'(\varepsilon)\) from van der
Corput Lemma), such that for most of the \( i\)'s in the interval
\( \{1,2,\ldots,I(\varepsilon,\mathbb{I}) \}\) (density of such \(
i \)'s should be at least \( 1 - \frac{\varepsilon}{3}\)) there
exists \( M(I(\varepsilon,\mathbb{I}),\mathbb{I},\varepsilon) \in
\mathbb{N} \), such that for every \( M \geq
M(I(\varepsilon,\mathbb{I}),\mathbb{I},\varepsilon) \)
\begin{equation}
\label{vdr_crp}
\left | \frac{1}{M} \sum_{m=1}^M
<v_m,v_{m+i}>_X\right| < \frac{\varepsilon}{2},
\end{equation}
for all \( \{i_1,\ldots,i_j\} \in S \).
\newline
 In our case we obtain
\[
\left | \frac{1}{M} \sum_{m=1}^M <v_m,v_{m+i}>_X\right|=
\]
\[
| \frac{1}{X} \sum_{x=1}^X \frac{1}{M} \sum_{m=1}^M
\psi^{1,j+1}_{\{r_1^1,c_1^1 i_1\},\ldots,\{r_j^1,c_j^1 i_j\},\{0,b_1
i \}}(a_1x+ b_1 m) \ldots \]
\[
\psi^{k,j+1}_{\{r_1^k,c_1^k i_1\},\ldots,\{r_j^k,c_j^k i_j\},\{0,b_k
i \}}(a_kx+ b_k m) |=\tilde{A}.
\]
Denote \( y = a_1x+b_1m \). Assume that \( (a_1,b_1) = d \). Denote
\[
\tilde{B}_{y,m} = \psi^{1,j+1}_{\{r_1^1,c_1^1
i_1\},\ldots,\{r_j^1,c_j^1 i_j\},\{0,b_1 i \}}(y)
 \ldots
\psi^{k,j+1}_{\{r_1^k,c_1^k i_1\},\ldots,\{r_j^k,c_j^k i_j\},\{0,b_k
i \}}(a_k'y + b_k'm),
\]
where \(
a_p' = \frac{a_p}{a_1} \), \(b_p' = b_p - a_p'b_1 \).
Now we rewrite \( \tilde{A} \) in the following way
\[ \tilde{A} =
\left|a_1 \frac{1}{Y} \left(\sum_{l=0}^{\frac{a_1}{d}-1}
\sum_{y\equiv dl \mod a_1}^Y \frac{1}{M} \sum_{m \equiv \phi(l) \mod
\frac{a_1}{d}}^M \tilde{B}_{y,m}
%
\right)\right|
 + \delta_{X,M}.
\]
Here \( \phi \) is the one to one function from
\(\mathbb{Z}_{\frac{a_1}{d}}\) onto itself, such that
 \( \phi(l)\frac{b_1}{d} \equiv l \mod \frac{a_1}{d} \)
for every \( 0 \leq l \leq \frac{a_1}{d}-1 \), \( Y = a_1 X \),
\(a_p',b_p'\) as above and \( \delta_{X,M} \) accounts for the fact
that for small \( y\)'s and \( y \)'s close to \( Y \) there is a
difference between elements that are taken in the expression for \(
\tilde{A} \) and in the expression on the right hand side of the
last equation. Nevertheless, we have \( \delta_{X,M} \rightarrow 0\)
if \( \frac{M}{X} \rightarrow 0 \).
\newline
It will suffice to prove (Cauchy-Schwartz inequality) that there exists
\( \mathbb{I}(\varepsilon,\delta) \in \mathbb{N} \), such that for
every \( \mathbb{I} \geq \mathbb{I}(\varepsilon,\delta) \) there
exists a set \( S \subset [-\mathbb{I},\mathbb{I}]^j \) of density
at least \( 1-\delta \) and there exist \(
I(\varepsilon,\mathbb{I}) \in \mathbb{N} \), \(
M(I(\varepsilon,\mathbb{I})) \in \mathbb{N} \), such that for
every \( M \geq M(I(\varepsilon,\mathbb{I}))\) there exists \(
X(M,\varepsilon) \in \mathbb{N} \) such that for every \( X \geq
X(M,\varepsilon) \), and for a set of \( i \)'s in the interval \(
\{1,2,\ldots,I(\varepsilon,\mathbb{I})\}\) of density \(1 -
\frac{\varepsilon}{3}\) we have
\begin{equation}
\label{sub_eq_lemma}
 a_1\frac{1}{Y}  \sum_{y\equiv dl \mod a_1}^Y
\left(\frac{1}{M}
 \sum_{m \equiv \phi(l) \mod
\frac{a_1}{d}}^M \tilde{C}_{y,m}
 \right)^2 < \left( \frac{\varepsilon d}{3a_1} \right)^2,
\end{equation}
for all \( 0 \leq l \leq \frac{a_1}{d}-1\), where
\[
\tilde{C}_{y,m} = \psi^{2,j+1}_{\{r_1^2,c_1^2
i_1\},\ldots,\{r_j^2,c_j^2 i_j\},\{0,b_2 i \}}(a_2'y + b_2'm)
 \ldots \]
 \[
\psi^{k,j+1}_{\{r_1^k,c_1^k i_1\},\ldots,\{r_j^k,c_j^k i_j\},\{0,b_k
i \}}(a_k'y + b_k'm).
\]
\newline
We rewrite the inequality (\ref{sub_eq_lemma}) for a fixed \( l \)
as follows:
\newline
Denote \( z \) and \( n \), such that \( y = z a_1 + dl \) and \( m
= n \frac{a_1}{d}+ \phi(l) \). As a result we obtain
\[
\frac{1}{Z} \sum_{z=1}^Z \left(\frac{d}{Na_1} \sum_{n=1}^N
\psi^{2,j+1}_{sh_2}\left(t_{n,z,l}^2\right)
 \ldots
\psi^{k,j+1}_{sh_k}\left(t_{n,z,l}^k\right)
  \right)^2=
\]
\[
\frac{1}{Z} \sum_{z=1}^Z \left(\frac{d}{Na_1} \sum_{n=1}^N
\psi^{2,j+1}_{sh_2}\left(a_2 z +c_2n + r_2\right)
 \ldots
\psi^{k,j+1}_{sh_k}\left(a_k z +c_kn + r_k\right)
  \right)^2 \doteq \tilde{D},
\]
where \( sh_p = \{\{r_1^p,c_1^p i_1\},\ldots,\{r_j^p,c_j^p
i_j\},\{0,b_p i \}\} \), \newline \( t_{n,z,l}^p = \frac{a_p (a_1z
+dl) + (a_1b_p-a_pb_1)(\frac{a_1}{d}n+\phi(l)) }{a_1}\), \( r_p =
\frac{a_pl + (a_1 b_p-a_p b_1) \phi(l)}{a_1} \), \newline \( c_p =
\frac{a_1 b_p - a_p b_1}{d} \neq 0 \), \( Z = \frac{Y}{a_1} \) and
\( N = \frac{Md}{a_1} \). The expression \( \frac{a_pl + (a_1b_p-a_p
b_1) \phi(l)}{a_1} \in \mathbb{Z} \) (from the condition on the
function \( \phi \)).

\noindent    From the conditions of the lemma we obtain for every \(
p \neq q , \, \, p,q > 1 \)
\[
\det \left( \begin{array}{cc}
a_p & c_p \\
a_q & c_q  \\
\end{array}
\right) = \frac{a_1 \det \left(\begin{array}{cc}
a_p & b_p \\
a_q & b_q  \\
\end{array}
\right)}{d} \neq 0.
\]
\newline
Therefore, we have \( \tilde{D} \) can be rewritten
\[
\tilde{D} = \frac{1}{Z} \sum_{z=1}^Z (\frac{1}{Na_1} \sum_{n=1}^N
\psi^{2,j+1}_{\{r_1^2,c_1^2 i_1\},\ldots,\{r_j^2,c_j^2
i_j\},\{r_2,b_2 i\}}\left(a_2 z +c_2n \right)
 \ldots
 \]
 \[
\psi^{k,j+1}_{\{r_1^k,c_1^k i_1\},\ldots,\{r_j^k,c_j^k
i_j\},\{r_k,b_k i\}}\left(a_k z +c_kn \right)
 )^2.
\]
By the induction hypothesis there exists \(
\mathbb{I}(\varepsilon,\delta) \in \mathbb{N} \)  big enough, such
that for every \( \mathbb{I} \geq \mathbb{I}(\varepsilon,\delta)
\) there exist a subset \( S \subset
[-\mathbb{I},\mathbb{I}]^{j+1} \) of density at least \( 1 -
\delta^2 \) and \( N(\mathbb{I},\varepsilon) \in \mathbb{N} \),
such that for every \( N \geq N(\mathbb{I},\varepsilon) \) there
exists \(Z(N,\mathbb{I},\varepsilon) \in \mathbb{N} \), such that
for every \( Z \geq Z(N,\mathbb{I},\varepsilon)\)
\begin{equation}
\label{ineq_last} \tilde{D} < \left( \frac{\varepsilon}{3a_1}
\right)^2,
\end{equation}
 for all \( \{i_1,\ldots,i_j,i\} \in S\).
\newline
For every \( (i_1,\ldots,i_j) \in [-\mathbb{I},\mathbb{I}]^j \) we
denote by \( S_{i_1,\ldots,i_j} \) the following subset of \(
[-\mathbb{I},\mathbb{I}] \)
\[
S_{i_1,\ldots,i_j} = \{ i \in [-\mathbb{I},\mathbb{I}] \,\, | \,\,
(i_1,\ldots,i_j,i) \in S \}.
\]
Then there exists a set \( T \subset [-\mathbb{I},\mathbb{I}]^j \)
of density at least \( 1 - \delta \), such that for every \(
(i_1,\ldots,i_j) \in T \) the density of \( S_{i_1,\ldots,i_j} \) is
at least \( 1 - \delta \). Let \( \delta < \frac{\varepsilon}{7} \)
 and \( \mathbb{I}
> \max_{l}{(\max{(I'(\varepsilon),\mathbb{I}(\varepsilon,\delta))})}
\) (\( I'(\varepsilon)\) is taken from van der Corput lemma). By
taking \( N(\mathbb{I},\varepsilon,\delta) \),  follows
 from the inequality (\ref{ineq_last}) that there exists \(
 M(\mathbb{I},\varepsilon,\delta) \in \mathbb{N}\), such that for every
 \( M \geq M(\mathbb{I},\varepsilon,\delta)\)
 there exists \( X(M,\mathbb{I},\varepsilon,\delta) \in \mathbb{N} \), such that for
 every \( X \geq X(M,\mathbb{I},\varepsilon,\delta) \) the
 inequality (\ref{vdr_crp}) holds
for every fixed \( (i_1,\ldots,i_j) \in T \) for a set of \( i
\)'s within the interval \(\{1,\ldots,\mathbb{I}\}\) of density at
least \( 1 - \frac{\varepsilon}{3} \). The lemma follows from the van der Corput lemma.

\hspace{12cm} \qed
\end{proof}

\noindent
\begin{proof}\textit{\( \openrm \)of Proposition  \rm{\ref{main_prop}}\( \closerm \)}
\newline
Denote by \( v_m(x) \doteq \xi_1(a_1 x + b_1 m) \ldots \xi_k(a_k x +
b_k m) \). Then for every \( i \in \mathbb{N} \)
\[
\left| \frac{1}{M} \sum_{m=1}^M <v_m,v_{m+i}>_X \right| =\]
\[
\left| \frac{1}{X} \sum_{x=1}^X  \frac{1}{M} \sum_{m=1}^M
\psi_{\{0,b_1i\}}^{1,1}(a_1 x+b_1 m) \ldots
\psi_{\{0,b_ki\}}^{k,1}(a_k x+b_k m)\right| \doteq \tilde{A},
\]
where functions \( \psi^{p,j} \)'s are autocorrelation functions of
\( \xi_p \)'s of the length \( j \).
\newline
Again, as in the proof of the lemma \ref{autocorrelation2} we denote
\( y = a_1 x + b_1 m \). We proceed with the analysis of the
expression \( (\tilde A) \) and by the same technique which was used
in the proof of the lemma \ref{autocorrelation2} we conclude the
following:
\newline
In order to prove that \( \tilde{A} < \frac{\varepsilon}{2}\) for a
set of \( i \)'s within the appropriate interval \(
\{1,2,\ldots,I(\varepsilon) \}\) it is sufficient to prove that
there exists \( I(\varepsilon) \in \mathbb{N}\) big enough and \(
N(I(\varepsilon),\varepsilon) \in \mathbb{N} \), such that for every
\( N \geq N(I(\varepsilon),\varepsilon) \) there exists \(
Z(N,\varepsilon) \in \mathbb{N} \), such that for every \( Z \geq
Z(N,\varepsilon)\)
\[
\frac{1}{Z} \sum_{z=1}^Z \left(\frac{1}{Na_1} \sum_{n=1}^N
\psi^{2,1}_{\{r_2,b_2 i\}}\left(a_2 z +c_2n \right)
 \ldots
\psi^{k,1}_{\{r_k,b_k i\}}\left(a_k z +c_k n \right)
  \right)^2 < \left( \frac{\varepsilon}{3a_1} \right)^2,
\]
for  a set \( i \)'s within the interval \(
\{1,\ldots,I(\varepsilon)\}\) of density \(
1-\frac{\varepsilon}{3}\) .
\newline
The last statement follows from lemma \ref{autocorrelation2}.
The proposition follows from the van der Corput lemma.

\hspace{12cm} \qed
\end{proof}
\subsection{Probabilistic constructions of WM sets}\label{lin_useful_constructions}
The goal of this section is to prove the necessity of the conditions
of theorem \ref{main_thm_lin} and  the following proposition is the main tool
for this task.
\begin{prop}
\label{3_constr_prop} Let \( a,b \in \mathbb{N} \), \( c \in \mathbb{Z} \) such that \( a
\neq b \). Then there exists a WM set \( A \) such that within it
the equation
\begin{equation}
\label{lin_eq_gen}
 ax=by + c
\end{equation}
is unsolvable, i.e., for every \( (x,y) \in A^2 \) we have \( ax
\neq by + c \).
\end{prop}
\begin{remark}
\textnormal{The proposition is a particular case of theorem \ref{main_thm_lin}. It is a crucial ingredient in proving
the necessity direction of the theorem in general.}
\end{remark}
\begin{proof}
Let \( S \subset \mathbb{N} \). We construct from \( S \) a new
set \( A_S \), such that within it the equation \( ax = by + c \) is
unsolvable. Without loss of generality, suppose that \( a < b \).
\newline
Assume \( (a,b)=1 \) (the general case follows easily).
The
equation \( ax = by +c \) is solvable only if \( x \equiv
\phi(a,b,c) \mod b \), where \( \phi(a,b,c)\,:\, 0 \leq \phi(a,b,c) < b \) is determined
uniquely (if the equation has a solution at all, otherwise any WM
set will provide an example). Let us denote \( l_0 \doteq  \phi(a,b,c)
\). We define inductively a sequence \( \{l_i\} \subset \mathbb{N} \cup \{0\}\).
If a pair \( (x,y) \) is a solution of
the equation and \( y \in b^i\mathbb{N}+l_{i-1} \) then there exists \( l_{i} \in \{0,1,\ldots,b^{i+1}-1\} \) such that
\( x \in b^{i+1} \mathbb{N} + l_i \).
\newline
We define the sets \( H_i \doteq b^i \mathbb{N} + l_{i-1} \, ; \, i \in \mathbb{N} \). We prove that for every
\( i \in \mathbb{N} \, : \, H_{i+1} \subset H_i \). All elements of \( H_{i+1} \) are in the same class
modulo \( b^i \). So, if we show for some \( x \in H_{i+1} \) that \( x \equiv l_{i-1} \mod(b^{i-1}) \)
then we are done. For \( i = 1 \) we know that if
\( y \in \mathbb{N}\) then \( x \, : \, ax = by + c \) has to be in \( H_1 \). Therefore for \( x \in H_2\) such that
there exists \( y \in H_1 \) such that \( ax = by + c \) we have that \( x \in H_1 \). Therefore,
we have shown that \( H_2 \subset H_1 \). For \( i > 1 \) there exists \( x \in H_{i+1} \) such that
there exists \( y \in H_{i} \) with \( ax = by +c \). By induction \( H_{i} \subset H_{i-1} \). Therefore,
the latter \( y \) is in \( H_{i-1} \). Therefore, by construction of \( l_i \)'s we have that
\( x \in H_{i} \). Thus, by aforementioned remark, we established \( H_{i+1} \subset H_i \).
We define  the sets \( B_i; 0 \leq i < \infty \):
\[
  B_0= \mathbb{N} \setminus H_1,
  \]
  \[
  B_1=  H_1 \setminus  H_2
\]
\[
\ldots
\]
\[
  B_i =  H_i \setminus  H_{i+1}
\]
\[
  \ldots
\]
\newline
Clearly we have \( B_i \cap B_j = \emptyset \, \,, \forall
i \neq j \) and \( |\mathbb{N} \setminus (\cup_{i=0}^{\infty} B_i )| \leq 1 \).
 The latter is because
for every \( i \) the second element (in the increasing order) of \( H_i \) is
\( \geq b^i \), therefore if the latter set would contain \( 2 \) elements then the second element (in the increasing order)
is unbounded.
\newline
We define \( A_S = \bigcup_{i=0}^{\infty} A_i \), where
\( A_i \)'s are defined in the following manner:
 \[ A_0 \doteq S \cap B_0 , C_0 \doteq A_0^c \cap B_0\]
\[ D_1 \doteq B_1 \setminus \{x \in B_1\,|\, ax \in b B_0 + c\},
 A_1 \doteq  \left( B_1 \cap \{x\,|\, ax \in b C_0 + c\} \right) \cup \left( D_1 \cap S \right),\]
 \[
 C_1 \doteq A_1^c \cap B_1 \]
\[ \ldots \]
\[
D_i \doteq B_i \setminus \{x \in B_i\,|\, ax \in b B_{i-1} + c\}, A_i = \left( B_i \cap
\{x\,|\, ax \in b C_{i-1} + c\} \right) \cup \left(D_i \cap S\right), \]
\[
C_i \doteq A_i^c \cap B_i
\]
\[
\ldots
\]
Here  it is worthwhile to remark that for every \( i \, : A_i \subset B_i \,\, and \,\, B_i = A_i \cup C_i\).
Therefore \( A_S \subset \cup_{i=0}^{\infty} B_i \).
\newline
If for some \( i \) we have \( y \in A_i \subset B_i\)  then \( x \, : \, ax = by + c \) satisfies \( ax \in b A_{i} + c \) and
by the construction of \( B_i \)'s we know that \( x \in B_{i+1}\) or
\( x \in \mathbb{N} \setminus (\cup_{i=0}^{\infty} B_i)\). In the first case
\( x \not \in A_{i+1} \Rightarrow x \not \in A_S \).
In the second case \(x \not \in A_S \).
\newline
Thus in \( A_S \) the equation
(\ref{lin_eq_gen}) is unsolvable. Our main claim is the following.

\noindent \textit{For almost every subset \( S \) of \( \mathbb{N} \) the set
\( A_S \) is a normal set.}

\noindent By normality we mean that the infinite binary sequence \( 1_{A_S}
\in \{0,1\}^{\infty} \) is a normal binary sequence. The probability
measure on subsets of \( \mathbb{N} \) is the product on \( \{0,1\}^{\infty} \) of
probability measures \( (\frac{1}{2},\frac{1}{2})\).

\noindent
The tool for proving the claim is the following easy lemma (for a proof see appendix, lemma
 \ref{norm_lemma}).
\newline
\textit{A subset \( A \) of natural numbers is a normal set  \(
\Leftrightarrow \) for any \(k \in (\mathbb{N} \cup \{0\}) \) and
any \( i_1 < i_2 <  \ldots < i_k \) we have
\begin{equation}
 \lim_{N \rightarrow \infty} \frac{1}{N} \sum_{n=1}^N \chi_A(n) \chi_A(n+i_1) \ldots \chi_A(n+i_k)
 =0,
 \label{eq:main_eq}
\end{equation}
where \( \chi_A(n) \doteq 2 1_A(n) - 1\).}
\newline
First of all, we denote by \( T_N = \frac{1}{N} \sum_{n=1}^N
\chi_{A_S}(n) \chi_{A_S}(n+i_1) \ldots \chi_{A_S}(n+i_k) \).
Because of randomness of \( S \), \(T_N \) is a random variable
(the probability on Borel subsets of \(\{0,1\}^{\mathbb{N}}\)). We will
prove that \( \sum_{N=1}^{\infty} E(T_{N^2}^2) < \infty \) and
this will imply that \( T_N \rightarrow_{N \rightarrow \infty} 0 \)
for almost every \( S \subset \mathbb{N} \).
\newline
\[
E(T_N^2) = \frac{1}{N^2} \sum_{n,m=1}^N E(\chi_{A_S}(n)
\chi_{A_S}(n+i_1) \ldots \chi_{A_S}(n+i_k) \chi_{A_S}(m)
\chi_{A_S}(m+i_1) \ldots \chi_{A_S}(m+i_k)).
\]
A unique possible element of complement of \( \cup_i B_i = \cap_{i=1}^{\infty} H_i\) doesn't effect
the normality of \(A_S \) and we assume without loss of generality that \( \cap_{i=1}^{\infty} H_i = \emptyset\),
thus \( \mathbb{N} = \cup_{i=0}^{\infty} B_i \).
For every number \( n \in \mathbb{N} \) we define the chain of \( n \), \( Ch(n) \), in the following way:
\newline
If \( n \in B_0 \), then \( Ch(n) = (n)\).
\newline
If \( n \in B_1 \), then two situations are possible. In the first one there exists a unique
\( y \in B_0 \) such that \( an = by +c \). We set \( Ch(n) = (n,y)=(n,Ch(y))\).
In the second situation we can not find such \( y \) from \( B_0\) and we set \(Ch(n) = (n) \).
\newline
If \( n \in B_{i+1}\), then again two situations are possible. In the first one there exists
\( y \in B_i \) such that \( an = by+c \). In this case we set
\(Ch(n)=(n,Ch(y))\). In the second situation there is no such \( y \) from no one of \(B_0,\ldots,B_i\).
In this case we set \( Ch(n) = (n) \). We define the length of \( Ch(n) \), \( l(n)\),
to be a number of elements
in \( Ch(n) \).

\noindent For every \( n \in \mathbb{N} \)  we define the ancestor of \( n \), \( a(n) \), to be the last
element of the chain of \( n \) (\(Ch(n)\)).
To determine  whether or not \( n \in A_S \) will
depend on whether \( a(n) \in S \). The exact relationship will depend on
the \( i \) for which \( n \in B_i \) and \( j \) for which \( a(n) \in B_j \) or in other words on length
of \( Ch(n)\):
\( \chi_{A_S}(n) = (-1)^{i-j} \chi_S(a(n)) = (-1)^{l(n)-1} \chi_S(a(n))\) (as proven below).


\noindent We say that  \( n \) is a descendant of \( a(n) \).

 \noindent We prove the formula  \( \chi_{A_S}(n) = (-1)^{i-j} \chi_S(a(n)) \), where \( i \) and  \(j \)
 are defined by \( n \in B_i \) and \( a(n) \in B_j \).
 \newline
 If \( i = 0 \Rightarrow j = 0\) and the formula is obvious.
 \newline
 For \( i > 0 \):  If \( j =i \) then the formula again is obvious. If \( j = i -1 \) then in
 case \( a(n) \in A_{i-1} \) we get that \( n \not \in A_i \) and in  case \( a(n) \not \in A_{i-1} \)
 we get that \( n \in A_i \). Therefore, we get \( \chi_{A_S}(n) = -\chi_{S}(a(n)) \). For general \( j < i -1 \)
 the argument is the same.

 \noindent It is evident that \( E(\chi_{A_S}(n_1) \ldots
\chi_{A_S}(n_k)) \neq 0 \Leftrightarrow \) every number \( a(n_i)
\) occurs an even number of times among numbers \( a(n_1),a(n_2),
\ldots, a(n_k)\).
\newline
We will bound the number of \( n,m\)'s inside the square \([1,N]
\times [1,N]\) such that \( E(\chi_{A_S}(n) \chi_{A_S}(n+i_1)
\ldots \chi_{A_S}(n+i_k) \chi_{A_S}(m) \chi_{A_S}(m+i_1) \ldots
\chi_{A_S}(m+i_k)) \neq 0 \).
\newline
For a given \( n \in [1,N] \) we will count all \( m \)'s inside
\( [1,N] \) such that for the ancestor of \( n \) there will be a chance
to have a twin among the ancestors of all \( n+i_1, \ldots,
n+i_k,m,m+i_1,\ldots,m+i_k \).
\newline
First of all it is obvious that in the interval \( [1,N]\) for a
given ancestor there can be at most \( \log_{\frac{b}{a}} N + C_1\) descendants, where \( C_1 \) is a
constant. For
all but a finite number of \( n \)'s it is impossible that among \( n+i_1,\ldots,n+i_k\)
there is the same ancestor as for \( n \). Therefore we should
focus on ancestors of the set \( \{ m, m+ i_1, \ldots, m+i_k\}\).
For a given \( n \) we might have at most \( (k+1) (\log_{\frac{b}{a}} N +C_1)\)
options for the number \( m \) to provide that for one of elements
of the set \( \{ m, m+ i_1, \ldots, m+i_k\}\) has the same
ancestor as \( n \). Therefore for most of \( n \in [1,N] \)
(except maybe a bounded number \( C_2  \) of \( n \)'s which depends
only  on \( \{ i_1,\ldots,i_k \}\) and doesn't depend on \( N \))
we have at most \( (k+1) (\log_{\frac{b}{a}} N +C_1)\) possibilities for \( m \)'s
such that
\[
E(\chi_{A_S}(n) \chi_{A_S}(n+i_1) \ldots \chi_{A_S}(n+i_k)
\chi_{A_S}(m) \chi_{A_S}(m+i_1) \ldots \chi_{A_S}(m+i_k) \neq 0.
\]
Thus we have
\[
E(T_N^2) \leq
\frac{1}{N^2} \left( \sum_{n=1}^N (k+1) (\log_{\frac{b}{a}} N +C_1)+ C_2 N \right) =
\frac{1}{N}((k+1)\log_{\frac{b}{a}} N + C_3),
\]
where \( C_3 \) is a constant.
This means that
\[
\sum_{N=1}^{\infty} E(T_{N^2}^2) < \infty.
\]
Therefore \( T_{N^2} \rightarrow_{N \rightarrow \infty} 0 \) for
almost every \( S \subset \mathbb{N} \). By lemma \ref{tec_lemma} it
follows that \( T_N  \rightarrow_{N \rightarrow
\infty} 0 \) almost surely.
\newline
In the general case, where \( a,b \) are not relatively prime, if \( c \)  satisfies (\ref{lin_eq_gen}) then it should be
 divisible by \( (a,b) \). Therefore by dividing the equation (\ref{lin_eq_gen}) by \( (a,b) \)
 we reduce the problem to the previous case.

\hspace{12cm} \qed
\end{proof}

\noindent
We will use the following notation:
\newline
Let W be a subset of \( \mathbb{Q}^n\). Then for any subset \(I=\{i_1,\ldots,i_p\} \subset \{1,2,\ldots,n\}\) we define
\[
Proj_{I} W = W_{I}= \{ (w_{i_1},\ldots,w_{i_p} \, | \, \forall w=(w_1,w_2,\ldots,w_n) \in W \}.
\]
The next step involves an  algebraic statement with a topological proof  which we have to
establish.
\begin{lemma}
\label{3_algebraic_lemma}
Let \( W \) be a non-trivial cone in \(
\mathbb{Q}^n \) which has the property that for every two vectors \(
\vec{x_1} = \{a_1,a_2,\ldots,a_n\}^t , \vec{x_2} =
\{b_1,b_2,\ldots,b_n\}^t \in W \)  there exist two coordinates \( 1
\leq i < j \leq n \) (depend on the choice of \(
\vec{x_1},\vec{x_2}\)) such that
\[
\det \left( \begin{array}{cc}
a_i & b_i \\
a_j & b_j  \\
\end{array}
\right) = 0,
\]
then there exist at least two coordinates \( i < j \) such that the
projection of \( W \) on these two coordinates is of dimension \(
\leq 1\) (\( dim_{\mathbb{Q}} \, span \, Proj_{i,j} W \leq 1\)).
\end{lemma}
\begin{proof}
First of all \( W \) has a non-empty interior in the topological space
\( V = Span W \). Assume that there no exist \( i \neq j \) such that the projection of \( W \) on these two coordinates is of dimension \( \leq 1 \). Without loss of generality we assume that \( V = \mathbb{Q}^n \).
Let fix an arbitrary non-zero element \( \vec{x} \in W \).  For every \( i,j: \, 1 \leq i < j \leq n \) we define
the subspace \( U_{i,j} = \{ \vec{v} \in V \, | \, Proj_{i,j} \vec{v} \in Span Proj_{i,j} \vec{x}\}\).
\newline
From assumptions of the lemma it follows   
that \( W = \cup _{i,j;1 \leq  i < j \leq n} (W \cap U_{i,j})\). For every \(i \neq j \) we obviously have that the interior of \( U_{i,j} \) is empty set. We get a contradiction because a finite union of sets with empty interior can not be equal to a set with non-empty interior.


\hspace{12cm} \qed
\end{proof}

\noindent
\begin{proof}\textit{(of theorem \ref{main_thm_lin}, \(\Rrightarrow\))}
\newline
First of all, we shift the affine space of solutions of equation
(\ref{lin_eq}) to obtain a vector subspace, denote it \( U \). The
linear space \( U \) must contain vectors with all positive
coordinates. Otherwise, the solution space can have only finitely
many positive solutions. Take any WM set and delete a finite number
of its elements we obtain a set in which the system is not solvable.
But removing a finite number of elements from a  WM set does not
affect the statistics of the remaining set; therefore, it will be
still a WM set. Thus, we can generate a WM set \( A \) in which the
equation (\ref{lin_eq}) is not solvable. The latter contradicts the
assumption that the system is solvable within every WM set.

\noindent  Denote by \( W = \{
\vec{v} \in U \, | \, \langle\vec{v},\vec{e_i}\rangle > 0 \, , \,
\forall \, i: \, 1 \leq i \leq k \} \). \( W \) is a non-trivial cone. By excluding all coordinates \( i : \, 1 \leq i \leq n \) for which we have \(Proj_i W =\{0\}\) we can assume that for every \(i: \, 1 \leq i \leq n \) we have 
\( Proj_i W \neq \{0\}\). By lemma \ref{3_algebraic_lemma} we deduce that 
there exist maximal subsets of coordinates \( F_1,\ldots,F_l \) such that for every \( p \in
\{1,2\ldots,l\}\) we have \( \forall i,j \in F_p \) the space \(
V_{i,j} \doteq Span W_{i,j} \) is one dimensional.






 \noindent We fix \( p \, : \, 1 \leq p \leq l \). We should show that the projection on \( F_p \) of all solutions of
(\ref{lin_eq}) is on a shifted diagonal, where a shift is the same
for all coordinates in \( F_p \). If the projection of \( W \) on
coordinates from \( F_p \) is not on a diagonal then there exist two
coordinates \( i < j \) from \( F_p \) such that \( W_{i,j} =
\{(ax,bx) \, | \, x \in \mathbb{N}\}\) for some \( a \neq b \)
natural numbers. Therefore the projection of the solutions space of
(\ref{lin_eq}) on \(i,j\) has the form \(\{(ax+f_1,bx+f_2) \, | \, x
\in \mathbb{N}\} \), where \( f_1,f_2 \) are integers. From
proposition \ref{3_constr_prop} it follows that for any \( a,b,c \),
where \( a \neq b \), there exists a WM set \( A \) such that the
equation \( ax = by + c \) is not solvable inside \( A \). This
proves the existence of a WM set \( A \) such that for every \( x
\in \mathbb{Z} \) we have \( (ax+f_1,bx+f_2) \not \in A^2 \) (we
take a WM set \( A \) such that the equation \( a x = by + (a f_2 -
b f_1) \) is unsolvable inside \( A \)). To prove that a shift is
the same quantity for all coordinates in \( F_p \) we merely should
know that for any natural number \( c \) there exists a WM set \(
A_c \) such that inside \( A_c \) the equation \( x - y = c \) is
not solvable. The last statement is easy to verify.


\noindent Denote by \( E = \{1,2,\ldots,k\} \backslash (F_1 \cup F_2 \cup \ldots \cup F_l) \). Then
 there exist two vectors \( \vec{x_1},\vec{x_2} \in W \) such that the projection of them on two arbitrary
   coordinates from \( E \) is two dimensional. Therefore the condition a) of the theorem
  \ref{main_thm_lin}
holds. Moreover, by the same argument that was used to extract
maximal subsets of coordinates \(F_1,\ldots,F_l\) and by preceding
remarks there exist  \( \vec{x_1},\vec{x_2}\) which satisfy
condition a) and additionally satisfy conditions b) and c) of the
theorem. This completes the proof.

\hspace{12cm} \qed

\end{proof}

\subsection{ Comparison with Rado's Theorem}\label{sub_sect_Rado}

We recall that the problem of solvability of a system of linear
equations for any finite partition of \( \mathbb{N} \) was solved
by Rado in  \cite{rado}. Such systems of linear equations
are called regular (or partition regular). Before citing
the theorem we would say that we may  expect regular
systems to be solvable as well inside every WM set. This is in fact the case and could be
shown directly, without use of theorem \ref{main_thm}, by the
technique of Furstenberg and Weiss
that was developed in their dynamical proof of Rado's theorem (see
\cite{furst3}). Instead of doing so, we obtain this result by use of
theorem \ref{main_thm_lin}.
\newline
First of all we should describe Rado's regular systems. We will need
a definition of the following object.
\begin{definition}
A rational \( p \times q \) matrix \( (a_{ij}) \) is said to be of
level \( l \) if the index set \( \{1,2,\ldots,q\}\) can be divided
into \( l \) disjoint subsets \( I_1,I_2,\ldots,I_l \) and rational
numbers \( c_j^r \) may be found for \( 1 \leq r \leq l \) and \( 1
\leq j \leq q \) such that the following relationships are
satisfied:
\[
\sum_{j \in I_1} a_{ij} = 0
\]
\[
\sum_{j \in I_2} a_{ij} = \sum_{j \in I_1} c_j^1 a_{ij}
\]
\[
\ldots
\]
\[
\sum_{j \in I_l} a_{ij} = \sum_{j \in I_1 \cup I_2 \cup \ldots \cup
I_{l-1}} c_j^{l-1} a_{ij}
\]
for \( i=1,2,\ldots,p\).
\end{definition}
\begin{theorem}\textit{(Rado)}
A system of linear equations is regular if and only if  for some \(
l \) the matrix \( (a_{ij})\) is of level \( l \) and it is
homogeneous, i.e. a system is of the form
\[
\sum_{j=1}^q a_{ij}x_j = 0, \hspace{0.5 in} i=1,2,\ldots,p.
\]
\end{theorem}
After recalling Rado's result we
are ready to demonstrate the following.
\begin{prop}
\label{Rado_analog} A regular system is solvable in every WM set.
\end{prop}
\begin{proof}
Let a system \( \sum_{j=1}^q a_{ij}x_j = 0, i=1,2,\ldots,p\) be
regular. We will use the fact that the system is solvable for any
finite partition of \( \mathbb{N} \). First of all, the set of
solutions of a regular system is a subspace of \( \mathbb{Q}^q \),
let us denote it \( V \). It is obvious that \( V \) contains
vectors with all positive components. If for some \( 1 \leq i < j
\leq q \) we have \( Proj_{i,j}^{+} V \) (where \( Proj_{i,j}^{+} V
= \{(x,y)|x,y \geq 0 \hspace{0.1 in} \& \hspace{0.1 in} \exists
\vec{v} \in V: \, <\vec{v},\vec{e_i}>=x \, , \,
<\vec{v},\vec{e_j}>=y \} \)) is contained in a line, then \(
Proj_{i,j}^{+} V \) is diagonal, i.e. is contained in \( \{(x,x)| x
\in \mathbb{Q}\} \). Otherwise, we can generate a partition of \(
\mathbb{N} \) into two disjoint sets \( S_1,S_2 \)  such that no \(
S_1^q \) and no \( S_2^q \) intersects \( V \):

\noindent  This partition is constructed by an iterative process.
Without loss of generality we may assume that the line is \( x = n y
\), where \( n \in \mathbb{N} \). The general case is treated in the
same way. We start with \( S_1 = S_2 = \emptyset \). Let \( 1 \in
S_1 \).
\newline
Then we "color" the infinite geometric progression \( \{ n^m \, | \, m \in \mathbb{N} \} \) (adding elements
 to either \( S_1\) or
\( S_2 \)) in such way that there is no
\( (x,y) \) on the line from \( S_1^2,S_2^2 \).
Then we take a minimal element from \( \mathbb{N} \) which is still uncolored. Call it \( a \). Then
we add \( a \) to \( S_1 \). And again we "color" \( \{ a n^m \, | \, m \in \mathbb{N} \} \).
\newline
By induction in this way we obtain a desired partition of \( \mathbb{N} \).

\noindent This contradicts the assumption that the given system is
regular.
\newline
If for all \( 1 \leq i < j \leq k \) we have \( \dim_{\mathbb{Q}}
span(Proj_{i,j}^{+} V) = 2 \) then by lemma \ref{3_algebraic_lemma}
it follows that there exist two vectors \( \vec{x_1},\vec{x_2} \in V
\) which satisfy all requirements of theorem \ref{main_thm_lin}.
Thus, in this case the system is solvable in every WM set.
\newline
Otherwise, let \( F_1,\ldots,F_l \) denote maximal subsets of indices such that for every \( p \in \{1,\ldots,l\} \)
we have for every \( i \neq j \, , \, i,j \in F_p: \, \dim_{\mathbb{Q}} span(Proj_{i,j}^{+} V) = 1 \).
Let \( E = \{1,2,\ldots,k\} \setminus (F_1 \cup \ldots \cup F_l) \). For every \( p: 1 \leq p \leq l \)  we choose
arbitrarily one representative index within \( F_p \) and denote it by \( j_p \) (\(j_p \in F_p \)). Then by passing to the subset
of indices \( I \doteq E \cup \{j_1,\ldots,j_l\}\) we can show by use lemma \ref{3_algebraic_lemma} that there
exist \( \vec{x_1},\vec{x_2} \in V\) with all positive coordinates such that for every
\( i \neq j \, , \, i,j \in I \) we have \( \dim_{\mathbb{Q}} Proj_{i,j} (span(\vec{x_1},\vec{x_2})) = 2 \).
The latter ensures that the vectors \( \vec{x_1},\vec{x_2}\) satisfy all requirements of theorem
\ref{main_thm_lin} and, therefore, the system is solvable in every WM set.

\hspace{12cm} \qed
\end{proof}
\newpage
\section
{An additive analog of polynomial multiple recurrence for WM Sets}\label{additive_analog_section}

\noindent We recall the notation which was introduced earlier.
\newline
\textbf{Notation:} \textit{The Hilbert space \( L^2(N) \) is the space of all real-valued
functions on the finite set \( \{1,2,\ldots,N\}\) endowed with the
following scalar product:
\[
<u,v>_{N} = \frac{1}{N} \sum_{n=1}^N u(n)v(n).
\]
We denote by \( \left\Vert u \right\Vert_{N} = \sqrt{<u,u>_{N}}\).}

\noindent
The following definition will be used extensively.
\begin{definition}
Polynomials \( p_1,\ldots,p_k \in \mathbb{Z}[n] \) are called
\textbf{essentially distinct} if the difference of every two of them
is non-constant polynomial.
\end{definition}

\noindent The main result of this chapter is the following
\newline
\textbf{Theorem \ref{main_thm}:} \textit{For every \( k \in \mathbb{N} \) the system
\begin{equation}
\label{add_eq_sect_add_analog}
  \left\{ \begin{array}{llll} x+y_1=p_1(z) \\
x+y_2 = p_2(z) \\
\ldots \\
x+y_k=p_k(z)
\end{array} \right.
\end{equation}
 is solvable within every WM set if 
\( \deg (p_1) = \deg(p_2) = \ldots = \deg (p_k) \), \( p_1,\ldots,p_k \) are essentially distinct and have positive
leading coefficients.}

\subsection{Orthogonality of polynomial shifts}\label{subsect_orthogonality_polynomial_shifts}


The following lemma is essentially the main tool in the proof of the foregoing theorem.
It is inspired by the analogous proposition \ref{main_prop} in section \ref{lin_proof_suff}.
\begin{lemma}
\label{first_orth_lemma}
Let \( A \subset \mathbb{N} \) be a  WM set   and assume that \( p_1,\ldots,p_k \in \mathbb{Z}[n] \) are essentially distinct
 polynomials with positive leading coefficients.  We set \( \xi(n)
= 1_{A}(n) - d(A)  \) for non-negative \( n \) and zero for \( n \leq 0 \), and we assume
 \( q(n) \in \mathbb{Z}[n] \) with a
positive leading coefficient, \( \deg( q) \geq \max_{1 \leq i \leq k}\deg (p_i) \) and for every
\( i  \, : \, 1 \leq i \leq k \) such that \( \deg(p_i) = \deg(q) \) we have
 that the leading coefficient
of \( q(n) \) is bigger than that of \( p_i \).
Then for  every \( \varepsilon > 0 \) there exists \( J(\varepsilon) \)
such that for every \( J \geq J(\varepsilon) \) there exists \(
N(J,\varepsilon) \) such that for every \( N \geq N(J,\varepsilon)
\) we have
\[
\left\Vert \frac{1}{J} \sum_{j=1}^J a_{N+j}\xi(n - p_1(N+j)) \xi(n -
p_2(N+j)) \ldots \xi(n - p_{k}(N+j)) \right\Vert_{q(N)} <
\varepsilon
\]
for every \( \{a_n\} \in \{0,1\}^{\mathbb{N}} \).
\end{lemma}
\begin{proof}
We prove this statement by using an analog of Bergelson's PET
induction, see \cite{berg_pet}. Let \( F = \{p_1,\ldots,p_k\}\) be a
finite set of polynomials and assume that the largest of the degrees
of \( p_i \) equals \( d \). For every \( i \, : \, 1 \leq i \leq d
\) we denote by \( n_i \) the number of different groups of
polynomials of degree \( i \), where two polynomials \(
p_{j_1},p_{j_2} \) of  degree \( i \) are in the same group if and
only if they have the same leading coefficient. We will say that \(
(n_1,\ldots,n_d)\) is the \textit{characteristic vector} of \( F \).

\noindent We prove a more general statement than the statement of the
lemma.

 \noindent Let \( \mathcal{F}(n_1,\dots,n_d) \) be the
family of all finite sets of essentially distinct polynomials having
characteristic vector \( (n_1,\ldots,n_d)\). Consider the following
two statements:
\newline
\( L(k ; n_1,\ldots,n_d) \): 'For every \(
\{g_1,\ldots,g_{n_1},q_1,\ldots,q_l\} \in
\mathcal{F}(n_1,\ldots,n_d)\), where \( d \leq \deg(q) \), \( q \) is increasing faster
than any \( q_i, \, i: \, 1 \leq i \leq l \) (the exact statement is formulated in lemma) and  \(
g_1,\ldots,g_{n_1} \) are linear polynomials, and every \(
\varepsilon,\delta
> 0 \) there exists \( H(\delta,\varepsilon) \in \mathbb{N} \) such
that for every \( H \geq  H(\delta,\varepsilon) \) there exists \(
J(H,\varepsilon) \in \mathbb{N} \) such that for every \( J
\geq J(H,\varepsilon) \) there exists \(
N(J,H,\varepsilon) \in \mathbb{N}\) such that for every \( N
\geq N(J,H,\varepsilon) \) for a set of \(
\{h_1\ldots,h_k\} \in [1 \ldots H]^k \) of density at least \( 1 -
\delta \) we have
\[
\Vert \frac{1}{J} \sum_{j=1}^J a_{N+j} \prod_{i=1}^{n_1}
\prod_{\epsilon \in \{0,1\}^k} \xi(n - g_i(N+j) - \epsilon_1 h_1 -
\ldots - \epsilon_k h_k) \prod_{i=1}^{l} \xi(n-q_i(N+j))
\Vert_{q(N)} < \varepsilon,
\]
for every \( \{a_n\} \in \{0,1\}^{\mathbb{N}} \)'.
\newline
\( L(k;\overline{n_1,\ldots,n_i},n_{i+1},\ldots,n_d)\): '\( L(k;n_1,\ldots,n_d) \) is valid for any \( n_1,\ldots,n_i \)'.

\noindent  Lemma \ref{first_orth_lemma} is the special case \(
L(0; \overline{n_1,\ldots,n_d})\), where  \( d \leq \deg(q) \) and the polynomial
\( q \) is increasing faster than all polynomials in the given family of polynomials
which has the characteristic vector \((n_1,\ldots,n_d)\). In
order to prove the latter  it is enough to establish \( L(k;1) \,
, \, \forall k \in \mathbb{N} \cup \{0\}\), and to prove the
following implications:
\[
  S.1_d: \,  L(k+1;n_1,n_2,\ldots,n_d) \Rightarrow L(k;n_1+1,n_2,\ldots,n_d); \,\, k,n_1,\ldots,n_{d-1} \geq 0, n_d \geq 1, d \geq 1
\]
\[
  S.2_{d,i}: \, L(0;\overline{n_1,\ldots,n_{i-1}},n_i,\ldots,n_d) \Rightarrow L(k;\underbrace{0,\ldots,0}_{i-1 \, zeros},n_i+1,n_{i+1},\ldots,n_d);
\]
\[
\hspace{2in} k;n_1,\ldots,n_{d-1} \geq 0,n_d \geq 1, d \geq i >1
\]
\[
  S.3_d: \, L(k;\overline{n_1,\ldots,n_d}) \Rightarrow
L(k;\underbrace{0,\ldots,0}_{d \, zeros},1), \,\,\, k \geq 0 \, , \,
d \geq 1
\]

\noindent We start with a proof of statement \( S.2_{d,i} \). Suppose that \( F \) is a finite set of essentially distinct
polynomials and assume that the characteristic vector of \( F \) equals
\newline
\( (\underbrace{0,\ldots,0}_{i-1 zeros}, n_i+1,n_{i+1},\ldots,n_d)
\). Fix any of the \( n_i + 1 \) groups of polynomials of degree
\( i \) and denote its polynomials by \( g_1,\ldots,g_m \). Denote
the remaining polynomials in \( F \) by \( q_1,\ldots,q_l \).
Because there are no linear polynomials among the polynomials of
\( F \) , we have to show the following:

\noindent \textit{Let the family \( F \doteq
\{g_1,\ldots,g_m,q_1,\ldots,q_l\}\) of polynomials with the
characteristic vector \( (\underbrace{0,\ldots,0}_{i-1
zeros},n_i+1,n_{i+1},\ldots,n_d)\), where \(\{g_1,g_2,\ldots,g_m\}
\in \mathbb{Z}[n] \) is one of the groups of \( F\) of the degree
\( i, \, i > 1 \). Let \( A \) be a WM set and denote by \( \xi \)
the normalized WM-sequence, i.e., \( \xi(n) = 1_A(n) - d(A) \, ,
\, \forall n \in \mathbb{N}\). For every \( \varepsilon, \delta
> 0 \) there exists \( H(\varepsilon,\delta) \in \mathbb{N} \)
such that for every \( H \geq H(\varepsilon,\delta) \) there
exists \( J(\varepsilon,H) \) such that for every \( J \geq
J(\varepsilon,H) \) there exists \( N(J,\varepsilon,H) \) such
that for every \( N \geq N(J,\varepsilon,H) \) for a set of \(
(h_1,\ldots,h_k) \in \{1,\ldots,H\}^k \) of density which is at
least \( 1 - \delta \) we have
\[
\Vert \frac{1}{J} \sum_{j=1}^J a_{N+j}\prod_{\epsilon \in \{0,1\}^k} \xi(n - \epsilon_1 h_1 -\ldots-\epsilon_k h_k) \xi(n - g_1(N+j))  \ldots
\xi(n - g_{m}(N+j))
 \]
  \[ \xi(n - q_{1}(N+j))\ldots
\xi(n - q_{l}(N+j))\Vert_{q(N)} < \varepsilon,
\]
for every \( \{a_n\} \in \{0,1\}^{\mathbb{N}}\) and with the
condition \( \deg(q) \geq d \) and \( q \) is increasing faster than any
\( q_i, \, i: \, 1 \leq i \leq l \).}

\noindent  Denote by
\[
 u_j(n) \doteq  a_{N+j}\xi(n - g_1(N+j))  \ldots \xi(n - g_{m}(N+j))
\]
\[
 \xi(n - q_{1}(N+j))\ldots \xi(n - q_{l}(N+j)) ,
\]
\[
w(n) = \prod_{\epsilon \in \{0,1\}^k} \xi(n - \epsilon_1 h_1
-\ldots-\epsilon_k h_k),
\]
\[
 v_j(n) = w(n) u_j(n),
\]
\[
  \hspace{2in}  n=1,\ldots,q(N).
\]
The sequence \( w(n) \) is bounded by \( 1 \) and therefore to prove
that \( \Vert \frac{1}{J} \sum_{j=1}^J v_j\Vert_{q(N)} \) is small
it is sufficient to show that  \( \Vert \frac{1}{J} \sum_{j=1}^J u_j
\Vert_{q(N)} \) is small.

\noindent We apply the van der Corput lemma (see lemma
\ref{vdrCorput} in appendix):

\[
   \frac{1}{J} \sum_{j=1}^J <u_j,u_{j+h}>_{q(N)} =
\]
\[
   \frac{1}{q(N)} \sum_{n=1}^{q(N)} \frac{1}{J} \sum_{j=1}^J
a_{N+j}\xi(n - g_1(N+j))  \ldots \xi(n-g_m(N+j))
\]
\[
\xi(n- q_1(N+j)) \ldots \xi(n-q_l(N+j))
\]
\[
a_{N+j+h}\xi(n - g_1(N+j+h))  \ldots \xi(n - g_{m}(N+j+h))
\]
\[
\xi(n - q_1(N+j+h)) \ldots \xi(n-q_l(N+j+h)) =
\]
\[
    \frac{1}{q(N)-g_1(N)} \sum_{n=1}^{q(N)} \xi(n) \frac{1}{J} \sum_{j=1}^J a_{N+j}a_{N+j+h} \xi(n - (g_2(N+j) - g_1(N+j)))
    \ldots
\]
\[
    \xi(n- ( g_m(N+j)-g_1(N+j) )) \xi(n- (q_1(N+j) - g_1(N+j))) \ldots
\]
\[
    \xi(n-(q_l(N+j)-g_1(N+j)))  \xi(n - (g_1(N+j+h) - g_1(N+j)))  \ldots
\]
\[
    \xi(n - (g_{m}(N+j+h) - g_1(N+j))) \xi(n - (q_1(N+j+h) - g_1(N+j))) \ldots
\]
\[
\xi(n-(q_l(N+j+h) - g_1(N+j))) + \delta_{N,J}=
\]
\[
   \frac{1}{q(N)} \sum_{n=1}^{q(N)-g_1(N)} \xi(n) \frac{1}{J} \sum_{j=1}^J b_{N+j}
   \xi(n - r_1(N+j)) \ldots \xi(n - r_{m-1}(N+j)) \xi(n - r_m(N+j)) \ldots
\]
\[
  \xi(n- r_{m+l-1}(N+j)) \xi(n - r_{m+l}(N+j)) \ldots \xi(n-r_{2m+l-1}(N+j))
\]
\[
  \xi(n-r_{2m+l}(N+j)) \ldots \xi(n - r_{2m+2l-1}(N+j)) + \delta_{N,J},
\]
where in the second equality we used a change of variable \( n
\leftarrow n=n-g_1(N+j) \), \( b_{N+j} = a_{N+j}a_{N+j+h} \), \(
\delta_{N,J} \rightarrow_{\frac{J}{N} \rightarrow 0} 0 \) and
\[
   \left\{ \begin{array}{llll} r_t(n) = g_{t+1}(n) - g_1(n) \, , \, t: \, 1 \leq t \leq m-1 \\
                               r_t(n) = q_{t-(m-1)}(n) - g_1(n) \, , \, t: \, m \leq t \leq m+l-1 \\
                               r_t(n) = g_{t-(m+l-1)}(n+h) - g_1(n) \, , \, t: \, m+l \leq t \leq 2m+l-1\\
                               r_t(n) = q_{t-(2m+l-1)}(n+h)-g_1(n) \, , \, t: \, 2m+l \leq t \leq 2m+2l-1.
\end{array} \right.
\]
For all but a finite number of \( h \)'s the polynomials \(
\{r_t(n)\}_{t=1}^{2m+2l-1} \) are essentially distinct, because \( i
> 1 \) and the polynomials \(g_1,\ldots,g_m,q_1,q_l\) are essentially
distinct. To see the last property we notice that if we take two
polynomials \( r_t\)'s from the same group (there are \( 4 \)
groups), then their difference is a
 non-constant
because the initial polynomials are essentially distinct. If we take
two polynomials from different groups then three cases are possible.
In the first case the difference of these polynomials is \( g_t(n+h)
-g_t(n) \) or \( q_t(n+h) - q_t(n) \) for some \( t\). We assume
that \( i > 1 \) therefore \( min_{1\leq t \leq l} \min( \deg (q_t),
\deg(g_1)) > 1 \) and from this it follows that  \( g_t(n+h) -g_t(n)
\) and \( q_t(n+h) - q_t(n) \) are non-constant polynomials. In the
second case we get for some \( t_1 \neq t_2 \): \( g_{t_1}(n+h) -
g_{t_2}(n) \) or \( q_{t_1}(n+h) - q_{t_2}(n) \). Here we note that
the map \( h \mapsto p(n+h) \) is an injective map from \(
\mathbb{N} \) to the set of essentially distinct polynomials, if \(
\deg(p) > 1 \). Thus, for all but a finite number of \( h \)'s we
get again a non-constant difference. In the third case we get for
some \( t_1,t_2 \): \(g_{t_1}(n+h) - q_{t_2}(n)\) or \( q_{t_1}(n+h)
- g_{t_2}(n)\). The resulting polynomial has the same degree as \(
q_t \).

\noindent  The characteristic vector of the set of polynomials \( \{r_1,\ldots,r_{2m+2l-1}\}\) has the form
\( (c_1,\ldots,c_{i-1},n_i,n_{i+1},\ldots,n_d)\). The polynomials from the second and the fourth group have the same
degree as \( q_t \) and the same leading coefficient as \( q_t \) if \( \deg(q_t) > \deg(g_1) \) and the leading
coefficient will be the difference of leading coefficients of \( q_t \) and \( g_1 \) if \( \deg(q_t) = \deg(g_1)\).
The polynomials from the first and the third group will be of degree smaller than \( \deg(g_1) \).

\noindent Applying \(
L(0;\overline{n_1,\ldots,n_{i-1}},n_i,\ldots,n_d) \) with the new polynomial
\( q(n)-g_1(n) \) which is increasing faster than all the polynomials
 \( \{r_t(n)\}_{t=1}^{2m+2l-1} \) and
the Cauchy-Schwartz inequality we get that for all but a finite number
of \( h \)'s and for every \( \varepsilon > 0 \) there exists \(
J(\varepsilon,h) \) such that for every \( J \geq J(\varepsilon,h)
\) there exists \( N(J,\varepsilon,h) \) such that for every \( N
\geq N(J,\varepsilon,h) \) we have
\[
  \left| \frac{1}{J} \sum_{j=1}^J <u_j,u_{j+h}>_{q(N)} \right| < \varepsilon,
\]
for every \( \{a_n\} \in \{0,1\}^{\mathbb{N}}\).

\noindent By the van der Corput lemma it follows that for every \(
\varepsilon > 0 \) there exists \( J(\varepsilon) \) such that for
every \( J \geq J(\varepsilon) \) there exists \( N(J,\varepsilon)
\) such that for every \( N \geq N(J,\varepsilon) \) we have
\[
  \left \Vert \frac{1}{J} \sum_{j=1}^J u_j \right \Vert_{q(N)} < \varepsilon,
\]
for every  \( \{a_n\} \in \{0,1\}^{\mathbb{N}}\). Thus we have shown
the validity of \( L(k;\underbrace{0,\ldots,0}_{i-1 zeros},
n_i+1,n_{i+1},\ldots,n_d) \).

\noindent We proceed with a proof of \( S.1_d \). We fix  the
 \( n_1 + 1 \) groups of the polynomials of degree \( 1 \) and denote its polynomials by
 \( g_1(n)=c_1n+d_1, \ldots, g_{n_1+1} = c_{n_1+1}n+d_1\).
 (By the assumption that
 all given polynomials are essentially distinct we get that in any group of degree \( 1 \) there is only one polynomial).
  The remaining  polynomials we denote by \( q_1,\ldots,q_l\). The set of polynomials
  \( \{g_1,\ldots,g_{n_1+1},q_1,\ldots,q_l\}\) has the characteristic vector \( ( n_1+1,n_2,\ldots,n_d)\).
  Again we  apply the van der Corput lemma. Let \( u_j(n) \)
  be defined as following
\[
 u_j(n) \doteq  a_{N+j}\prod_{i=1}^{n_1+1} \prod_{\epsilon \in \{0,1\}^k}\xi(n - g_i(N+j) - \epsilon_1 h_1 - \ldots - \epsilon_k h_k)
 \prod_{i=1}^l \xi(n - q_i(N+j)),
\]
\[
  \hspace{2in}  n=1,\ldots,q(N).
\]
Then we have
\[
  \frac{1}{J} \sum_{j=1}^J <u_j,u_{j+h}>_{q(N)} =
\]
\[
  \frac{1}{q(N)} \sum_{n=1}^{q(N)} \frac{1}{J} \sum_{j=1}^J a_{N+j} a_{N+j+h}
\]
\[
  \prod_{i=1}^{n_1+1} \prod_{\epsilon \in \{0,1\}^k}\xi(n - g_i(N+j) - \epsilon_1 h_1 - \ldots - \epsilon_k h_k)
 \prod_{i=1}^l \xi(n - q_i(N+j))
 \]
\[
 \prod_{i=1}^{n_1+1} \prod_{\epsilon \in \{0,1\}^k}\xi(n - g_i(N+j+h) - \epsilon_1 h_1 - \ldots - \epsilon_k h_k)
 \prod_{i=1}^l \xi(n - q_i(N+j+h))=
\]
\[
  \frac{1}{q(N)-g_1(N)} \sum_{n=1}^{q(N)} \prod_{\epsilon \in \{0,1\}^k} \xi(n - \epsilon_1 h_1 - \ldots - \epsilon_k h_k )
  \xi(n - \epsilon_1 h_1 - \ldots - \epsilon_k h_k - c_1 h)
\]
\[
  \frac{1}{J} \sum_{j=1}^J b_{N+j}
  \prod_{i=1}^{n_1} \prod_{\epsilon \in \{0,1\}^k} \xi(n-(c_{i+1}-c_1)(N+j) - (d_{i+1}-d_1) - \epsilon_1 h_1 - \ldots - \epsilon_k h_k)
\]
\[
\prod_{i=1}^{n_1} \prod_{\epsilon \in \{0,1\}^k}
\xi(n-(c_{i+1}-c_1)(N+j) - (d_{i+1}-d_1) - \epsilon_1 h_1 -
\ldots - \epsilon_k h_k - c_{i+1}h)
\]
\[
 \prod_{i=1}^l \xi(n-(q_i(N+j)-g_1(N+j))) \prod_{i=1}^l \xi(n-(q_i(N+j+h)-g_1(N+j))) +\delta_{N,J},
\]
where in the second equality we made a change of variable \( n
\leftarrow n - g_1(N+j) \) and  \( b_{N+j} = a_{N+j}a_{N+j+h} \, ,
\, \delta_{N,J} \rightarrow_{\frac{J}{N} \rightarrow 0} 0 \).

\noindent Denote by \( r_i(n) = (c_{i+1}-c_1)n + (d_{i+1}-d_1) \, , \, i: 1 \leq i \leq n_1 \),
\(  s_i(n) = q_i(n) - g_1(n) \, , \, t_i(n) = q_i(n+h) - g_1(n) \, , \, i: 1 \leq i \leq l \). Then the last expression may
be rewritten as
\[
  \frac{1}{q(N)} \sum_{n=1}^{q(N)-g_1(N)} \prod_{\epsilon \in \{0,1\}^k} \xi(n - \epsilon_1 h_1 - \ldots - \epsilon_k h_k)
  \xi(n - \epsilon_1 h_1 - \ldots - \epsilon_k h_k- c_1 h)
\]
\[
  \frac{1}{J} \sum_{j=1}^J b_{N+j} \prod_{i=1}^{n_1}  \prod_{\epsilon \in \{0,1\}^k}
  \xi(n-r_i(N+j) - \epsilon_1 h_1 - \ldots - \epsilon_k h_k)
\]
\[
  \xi(n-r_i(N+j) - \epsilon_1 h_1 - \ldots - \epsilon_k h_k -c_{i+1}h)
\]
\[
  \prod_{i=1}^l \xi(n - s_i(N+j)) \xi(n-t_i(N+j)) + \delta_{N,J}  \doteq E1 + \delta_{N,J}.
\]
For every \( i \, : \, 1 \leq i \leq l \) the polynomials \( s_i,t_i
\) are in the same group (have the same degree and the same leading
coefficient), therefore the characteristic vector of the family \(
\{s_1,t_1,\ldots,s_l,t_l\}\) is the same as of the family \(
\{s_1,s_2,\ldots,s_l\}\) and , obviously, the characteristic vector
of the latter family is the same as of the family \(
\{q_1,q_2,\ldots,q_l\}\) and is equal to \((0,n_2,n_3,\ldots,n_d)\).
Again the polynomial \( q(n)-g_1(n) \) is increasing faster than any polynomial in the family
\(\{s_1,t_1,\ldots,s_l,t_l\}\) .
By use of \( L(k+1;n_1,\ldots,n_d) \) and the Cauchy-Schwartz
inequality we show that \( |E1| \) is arbitrarily small for a set of
arbitrarily large density of \( (h_1,\ldots,h_k,h)\)'s. Therefore,
by the van der Corput lemma we deduce the  validity of \(
L(k;n_1+1,n_2,\ldots,n_d)\).

\noindent The proof of \( S.3_d \) goes exactly in the same way as that 
of \( S.2_{d,i} \).

\noindent Proof of \( L(k;1) \, , \, \forall k \in \mathbb{N} \cup
{0}\):

\noindent  Assume that \( g_1(n) = c_1 n + d_1 \, , \, c_1 > 0 \) and \( \deg (q) > 1 \).
We show that

\noindent \textit{For every \( \varepsilon, \delta > 0\) there exists
\( H(\delta,\varepsilon) \in \mathbb{N} \) such that
for every \( H \geq H(\delta,\varepsilon) \) there exists \( J(H,\varepsilon) \in \mathbb{N} \) such that for every
\( J \geq J(H,\varepsilon)\) there exists \( N(J,H,\varepsilon) \) such that for every \( N \geq N(J,H,\varepsilon) \)
we have for a set of \( (h_1,\ldots,h_k) \in \{1,\ldots,H\}^k \) of density which is at least \( 1- \delta \) the
following
\[
  \left \Vert \frac{1}{J} \sum_{j=1}^J a_{N+j} \prod_{\epsilon \in \{0,1\}^k}
  \xi(n - g_1(N+j) - \epsilon_1 h_1 - \ldots - \epsilon_k h_k)\right \Vert_{q(N)} < \varepsilon
\]
for every \( \{a_n \} \in \{0,1\}^{\mathbb{N}} \).
}

\noindent We recall that to a WM set \( A \) is associated the weakly-mixing system \( (X_{\xi},\mathbb{B},T,\mu)\).
We define
the function \( f \) on \( X_{\xi}\)  by the following rule:
 \( f(\omega) = \omega_0 \, , \, \omega = \{\omega_0,\ldots,\omega_n,\ldots\} \in X_{\xi}\). It is evident that \( f \)
is continuous and \( \int_{X_{\xi}} f(x) d\mu(x) = 0 \). By genericity of the point \( \xi \in X_{\xi} \) we get
\[
\frac{q(N)}{q(N)-g_1(N)}\left \Vert \frac{1}{J} \sum_{j=1}^J a_{N+j} \prod_{\epsilon \in \{0,1\}^k}
 \xi(n - g_1(N+j) - \epsilon_1 h_1 - \ldots - \epsilon_k h_k)\right \Vert_{q(N)}^2 \rightarrow_{N \rightarrow \infty}
\]
\begin{equation}
\label{wm_discr_cube}
  \int_{X_{\xi}}  \left( \frac{1}{J} \sum_{j=1}^J a_{N+J+1-j}  T^{c_1 j}
  \left( \prod_{\epsilon \in \{0,1\}^k} T^{\epsilon_1 h_1 +\ldots +\epsilon_k h_k} f(x) \right)\right)^2 d\mu(x).
\end{equation}
Denote by \( g_{h_1,\ldots,h_k} \) the following function on \( X_{\xi}\):
\[
  g_{h_1,\ldots,h_k}(x) = \prod_{\epsilon \in \{0,1\}^k} T^{\epsilon_1 h_1 + \ldots + \epsilon_k h_k} f(x) \, ,
  \, \forall x \in X_{\xi}.
\]
Then we use the following statement which can be viewed as a corollary of theorem \( 13.1 \) of Host and Kra in \cite{host-kra}
(\( \int_{X_{\xi}} f(x) d\mu(x) = 0 \)).

\noindent \textit{For every \( \varepsilon, \delta > 0 \) there exists \( H(\delta,\varepsilon) \in \mathbb{N} \)
 such that for every
\( H \geq H(\delta,\varepsilon) \) for a set of \( (h_1,\ldots,h_k) \in \{1,\ldots,H\}^k\) which has density at least
\( 1- \delta \) we have
\[
   \left| \int_{X_{\xi}} g_{h_1,\ldots,h_k}(x) d\mu(x) \right| < \varepsilon.
\]}

\noindent

\noindent Let \( \varepsilon, \delta  > 0 \). By the foregoing
statement there exists \( H(\delta,\varepsilon) \in \mathbb{N} \)
such that for every \( H \geq H(\delta,\varepsilon)  \) the set of
those \( (h_1,\ldots,h_k) \in \{1,\ldots,H\}^k \) such that
\[
  \left| \int_{X_{\xi}} g_{h_1,\ldots,h_k}(x) d\mu(x) \right| < \frac{\varepsilon}{4}
\]
has density at least \( 1 - \delta \).

\noindent Lemma \ref{wm_sub_lemma1} implies that there exists
\( J(H,\varepsilon) \in \mathbb{N} \) such that for every \( J \geq  J(H,\varepsilon)  \)
 we have
\[
  \left\Vert \frac{1}{J} \sum_{j=1}^J b_j T^{c_1 j} \left( g_{h_1,\ldots,h_k}(x) - \int_{X_{\xi}} g_{h_1,\ldots,h_k}(x) d\mu(x)\right) \right \Vert_{L^2(X_{\xi})} < \frac{\varepsilon}{4}
\]
for any sequence \( \{b_n\} \in \{0,1\}^{\mathbb{N}} \).

\noindent Therefore, by merging  the two last statements we conclude that
there exists \( H(\delta,\varepsilon) \in \mathbb{N} \) such that
for every \( H \geq H(\delta,\varepsilon) \)  there exists \(
J(H,\varepsilon) \in \mathbb{N} \) such that for every \( J \geq
J(H,\varepsilon) \)  and for a set of \( (h_1,\ldots,h_k) \in
\{1,\ldots,H\}^k \) which has density at least \( 1 - \delta \) we
have
\[
  \left\Vert \frac{1}{J} \sum_{j=1}^J b_j T^{c_1 j} g_{h_1,\ldots,h_k}(x)\right \Vert_{L^2(X_{\xi})} < \frac{\varepsilon}{2}
\]
for any sequence \( \{b_n\} \in \{0,1\}^{\mathbb{N}} \).

\noindent Finally, by use of (\ref{wm_discr_cube}), the fact that
\( \lim_{N \rightarrow \infty} \frac{q(N)}{q(N)-g_1(N)} > 0 \)  and the last
statement we deduce the validity of \( L(k;1) \).

 \hspace{12cm} \qed
\end{proof}

\noindent The next lemma is a simple consequence of the previous one
and is used in the next subsection to prove  theorem \ref{main_thm}.
\begin{lemma}
\label{main_lemma}
 Let \( A \subset \mathbb{N} \) be a WM set and \( p_1,\ldots,p_k \in \mathbb{Z}[n] \) are essentially distinct
 polynomials of the same degree \( d \) greater than \( 1 \), with positive leading coefficients such
  that \( p_1(n) > p_i(n) , \, \forall 1 < i \leq k\) for sufficiently large \( n \). Then
for every \( \varepsilon > 0 \) there exists \( J(\varepsilon) \)
such that for every \( J \geq J(\varepsilon) \) there exists \(
N(J,\varepsilon) \) such that for every \( N \geq N(J,\varepsilon)
\) we have
\[
\left\Vert \frac{1}{J} \sum_{j=1}^J a_{N+j}\xi( p_1(N+j) - n) \xi(
p_2(N+j) - n) \ldots \xi(p_{k}(N+j)-n) \right\Vert_{p_1(N)} <
\varepsilon
\]
for every \( \{a_n\} \in \{0,1\}^{\mathbb{N}} \), where \( \xi(n) =
1_A(n) - d(A) \) for non-negative \( n \)'s and zero for \( n \leq 0 \).
\end{lemma}
\begin{remark}
\textnormal{The lemma is true for the linear case as well, but a proof demands additional technical efforts.}
\end{remark}
\begin{proof}
For a family of polynomials \( F = \{p_1,\ldots,p_k\}\) with a
maximal degree \( d \) denote by \( n_d \) the number
of different leading coefficients of polynomials of degree \( d \) from the family \( F \).
%
%

\noindent  As in the proof of lemma \ref{first_orth_lemma} we fix one of the groups of polynomials of degree \( d \) (all polynomials in the same group have the same
 leading coefficient). Assume that the group \(\{ g_1,\ldots,g_m \}\) has the maximal leading coefficient among all polynomials
 \( p_1,\ldots,p_k \). The rest of the polynomials we denote by \( q_1,\ldots,q_l \). Without loss of generality assume that
 \( p_1 = g_1, \ldots, p_m = g_m \).
Denote by
\( u_j(n) \, \, , \, \, 1 \leq n \leq p_1(N) \) the following
expression
\[
u_j(n) = a_{N+j}\xi( p_1(N+j) - n) \xi( p_2(N+j) - n) \ldots
\xi(p_{k}(N+j)-n).
\]
For \( u_j \)'s  we get
 \[
\frac{1}{J} \sum_{j=1}^J <u_j,u_{j+h}>_{p_1(N)} = \frac{1}{p_1(N)}
\sum_{n=1}^{p_1(N)} \frac{1}{J} \sum_{j=1}^J a_{N+j}\xi( p_1(N+j) -
n) \ldots \]
\[
\xi(p_{k}(N+j)-n) a_{N+j+h}\xi( p_1(N+j+h) - n) \ldots
\xi(p_{k}(N+j+h)-n) =
\]
\[
 \frac{1}{p_1(N)}  \sum_{n=1}^{p_1(N)} \xi(n) \frac{1}{J} \sum_{j=1}^J b_{N+j} \prod_{i=1}^{m-1} \xi(n-( p_1(N+j) - p_{i+1}(N+j)))
\]
\[
 \prod_{i=1}^{l}\xi(n - ( p_1(N+j) - q_i(N+j))) \prod_{i=1}^m  \xi(n-(p_1(N+j) - p_{i}(N+j+h)))
 \]
 \[ \prod_{i=1}^{l}\xi(n - (p_1(N+j) - q_i(N+j+h)))
 +\delta_{J,N},
\]
where \( b_{n} = a_n a_{n+h} \) and \( \delta_{J,N} \rightarrow_{\frac{J}{N} \rightarrow 0} 0 \).

\noindent Denote by \( r_i(n) =  p_1(n) - q_i(n) \, ; \, s_i(n) =  p_1(n) - q_i(n+h)\, , \, i: 1 \leq i \leq l \) and
\( t_i(n)  = p_1(n) - p_i(n) \, ; \, f_i(n) = p_1(n) - p_i(n+h) \, , \, i: 1 \leq i \leq m \). Then for all
 but a finite number of \( h \)'s  the polynomials \( \tilde{F} \doteq \{r_1,\ldots,r_l,s_1,\ldots,s_l,t_2,\ldots,t_m,f_1,\ldots,f_m \} \)
 are essentially distinct and among them the polynomials of degree \( d \) have \( n_d \) different leading coefficients.
 Therefore by lemma \ref{first_orth_lemma} for all but a finite number of \( h \)'s
 the following expression is as small as we wish for appropriately chosen \( J,N \).
\[
  \Vert \frac{1}{J} \sum_{j=1}^J b_{N+j} \prod_{i=1}^{m-1} \xi(n-t_{i+1}(N+j))  \prod_{i=1}^{l}\xi(n - r_i(N+j))
  \]
  \[
  \prod_{i=1}^m  \xi(n-f_i(N+j)) \prod_{i=1}^{l}\xi(n - s_i(N+j)) \Vert_{p_1(N)}.
\]
Finally by Cauchy-Schwartz inequality and van der Corput's lemma
we get the desired conclusion. \hspace{12cm} \qed
\end{proof}

\subsection{Proof of theorem \ref{main_thm} 
}\label{subsect_proof_suff_theorem_main_thm}


\begin{proof} \textit{(
of  theorem \ref{main_thm})}
\newline
The linear case (the degree of the polynomials is \( 1 \)) follows
from theorem \ref{main_thm_lin}:

\noindent Denote \( p_i(z) = c_i z + d_i \, , \, \forall i: \, 1
\leq i \leq k \). We choose the following order of the variables
\((x,z,y_1,\ldots,y_k)\). The the affine space of the solutions of
the additive system (\ref{add_eq_sect_add_analog}) is \( \{(x,z,
x- c_1 z - d_1, \ldots, x - c_k z - d_k) \, | \, x,z \in
\mathbb{Q}\}\). If we take the vectors \(\vec{x_1} =
(1,0,1,1,\ldots,1), \vec{x_2} = (0,1,c_1,c_2,\ldots,c_k), \vec{f}
= (0,0,-d_1,-d_2,\ldots,-d_k) \) then all the requirements of
theorem \ref{main_thm_lin} regarding the system
(\ref{add_eq_sect_add_analog}) are valid. Thus the system
(\ref{add_eq_sect_add_analog}) is solvable within every WM set.

\noindent Assume we have an arbitrary WM set \( A \) and \( k \)
essentially distinct polynomials \( p_1,\ldots,p_k \in
\mathbb{Z}[n] \) (a difference of any two of them is a non
constant polynomial) of the same degree \( d > 1\) with  positive
leading coefficients and assume that for sufficiently large \( n
\)'s we have \( p_1(n) > p_i(n) \, , \, \forall i: \,  2 \leq i
\leq k \). Let's define the set \( F \) of all \( z \)'s where the
statement of the theorem fails, namely,
\[
F \risingdotseq \{ z \in \mathbb{N} \, | \, for \,\, any \,\,
(x,y_1,\ldots,y_k) \in A^{k+1} \,\, the \,\, system
\,\,\text{\ref{add_eq_sect_add_analog}} \,\, fails\,\, to \,\,
hold\}.
\]
We shall prove that \( d^{*}(F) = 0 \). Since \( d(A) > 0 \) we
can find \( z \in A, z \not \in F \) and this will yield a
solution to (\ref{add_eq_sect_add_analog}).

\noindent Denote by \( \{a_n\} \) the indicator sequence of \( F
\), i.e., \( a_n = 1_F(n) \). We define the sequence \( \xi \) to
be a normalized indicator sequence of \( A \): \( \xi(n) = 1_A(n)
- d(A) \, , \, n \in \mathbb{N} \) and zero for non-positive values of \( n \),
 where \( d(A) \) is the density of \( A \) which
exists.
\newline
We define the expression \(B_{N,J} \) to be
\begin{equation}
\label{def_B_N_J}
B_{N,J} \risingdotseq \frac{1}{p_1(N)} \sum_{n=1}^{p_1(N)}
\frac{1}{J} \sum_{j=1}^J a_{N+j} 1_A(n) 1_A(p_1(N+j) - n)
\end{equation}
\[
1_A(p_2(N+j) - n) \ldots 1_A(p_{k-1}(N+j)-n) \xi(p_k(N+j)-n).
\]

\noindent     Suppose that we have \( d^{*}(F) > 0\). Then there
exist intervals \( I_{l,J} = [u_{l,J}+1,u_{l,J} + J] \) (for \( J
\) big enough) such that \( u_{l,J} \rightarrow _{l \rightarrow
\infty} \infty \) and \( \frac{|F \cap I_{l,J}|}{J}
> \frac{d^{*}(F)}{2} \) for every \(l\) and \(J\) big enough.
By induction on \( k \) and \( i \) we  prove the validity of the following
claim.

\noindent \textbf{Claim 1:} \textit{For every \(i \, : \, 0 \leq i
\leq k-1\) and every \( \varepsilon > 0 \) there exist \(J,l \) big
enough such that
\[
|\frac{1}{p_1(u_{l,J})} \sum_{n=1}^{p_1(u_{l,J})} \frac{1}{J}
\sum_{j=1}^J b_{u_{l,J}+j} 1_A(n) 1_A(p_1(u_{l,J}+j) - n)\ldots
\]
\[
1_A(p_i(u_{l,J}+j) - n) \xi(p_{i+1}(u_{l,J}+j) - n)
 \ldots \xi(p_{k}(u_{l,J}+j)-n) | <
 \varepsilon
\]
for every \(\{0,1\}\)-valued sequence \(\{b_n\}\). }

\noindent A proof of claim \( 1 \) is by induction on \( i \) and \( k \).
\newline
In the sequel we use the notation \( <1_A,f(n)>_N \), where \( f(n) \) is defined for all
\( n=1,2,\ldots,N \); which has the same meaning
as \( <1_A,f>_N = \frac{1}{N} \sum_{n=1}^N 1_A(n) f(n)\).
\newline
For \( i = 0 \) and every \( k \) the statement is exactly of lemma \ref{main_lemma}. For every
\( i < k-1 \) we
will prove the statement of the claim for \( i +1 \) and \( k \) provided
the statement for \( i \) and \( k \), and for \( i \), \( k - 1\):
\[
|\frac{1}{p_1(u_{l,J})} \sum_{n=1}^{p_1(u_{l,J})} \frac{1}{J}
\sum_{j=1}^J b_{u_{l,J}+j} 1_A(n) 1_A(p_1(u_{l,J}+j) - n)\ldots
\]
\[
1_A(p_i(u_{l,J}+j) - n) 1_A(p_{i+1}(u_{l,J}+j) - n)
\xi(p_{i+2}(u_{l,J}+j) - n) \ldots \xi(p_{k}(u_{l,J}+j)-n)
| =
\]
\[
| <1_A,\frac{1}{J} \sum_{j=1}^J b_{u_{l,J}+j}
1_A(p_1(u_{l,J}+j) - n) \ldots
\]
\[1_A(p_i(u_{l,J}+j) - n)
(\xi(p_{i+1}(u_{l,J}+j) - n) + d(A)) \xi(p_{i+2}(u_{l,J}+j) - n)
\ldots
\]
\[\xi(p_{k}(u_{l,J}+j)-n)
>_{p_1(u_{l,J})}| \leq
\]
\[
| <1_A,\frac{1}{J} \sum_{j=1}^J b_{u_{l,J}+j}
1_A(p_1(u_{l,J}+j) - n) \ldots
\]
\[1_A(p_i(u_{l,J}+j) - n)
\xi(p_{i+1}(u_{l,J}+j) \xi(p_{i+2}(u_{l,J}+j) - n)
\ldots
\]
\[\xi(p_{k}(u_{l,J}+j)-n)
>_{p_1(u_{l,J})}| +
\]
\[
d(A) | <1_A,\frac{1}{J} \sum_{j=1}^J b_{u_{l,J}+j}
1_A(p_1(u_{l,J}+j) - n) \ldots
\]
\[1_A(p_i(u_{l,J}+j) - n)
 \xi(p_{i+2}(u_{l,J}+j) - n)
\ldots
\]
\[\xi(p_{k}(u_{l,J}+j)-n)
>_{p_1(u_{l,J})}| < \varepsilon,
\]
for big enough \(J,l\).
The first summand is small by the statement of the claim for \( i \) and \( k \), and
the second summand is small by the statement of the claim for \( i \) and \( k -1 \).
This ends the proof of claim \( 1 \).

\noindent  We will use the statement of claim \( 1 \) for \( i = k-1 \) and we call the statement claim \( 2 \).

\noindent \textbf{Claim 2:} \textit{For every \(\varepsilon > 0 \)
there exist \( J,l \) big enough such that  the expression
\[
|\frac{1}{p_1(u_{l,J})} \sum_{n=1}^{p_1(u_{l,J})} \frac{1}{J}
\sum_{j=1}^J b_{u_{l,J}+j} 1_A(n) 1_A(p_1(u_{l,J}+j) - n) \ldots\]
\[
  1_A(p_{k-1}(u_{l,J}+j)-n) \xi(p_k(u_{l,J}+j)-n) | <
 \varepsilon
\]
for every \{0,1\}-valued sequence \(\{b_n\}\). }

\noindent The next statement enables us to conclude about a boundedness away from zero
of \( B_{u_{l,J},J} \).

\noindent \textbf{Claim 3:} \textit{For every \( \delta > 0 \)
for big enough \( J,l \)  the expression
\[
\frac{1}{p_1(u_{l,J})} \sum_{n=1}^{p_1(u_{l,J})} \frac{1}{J}
\sum_{j=1}^J b_{u_{l,J}+j} 1_A(n) 1_A(p_1(u_{l,J}+j) - n)
 \ldots 1_A(p_{k}(u_{l,J}+j)-n)
\]
is bigger than \( c (1 - \delta) d^{k+1}(A) \frac{d^{*}(F)}{3} \),
where \( c = \min_{2 \leq i \leq k-1} \frac{c_i}{c_1} \) (\( c_i
\) is a leading coefficient of polynomial \( p_i \)) for every
\(\{0,1\}\)-valued sequence \(\{b_n\}\) which has density bigger than
\(\frac{d^{*}(F)}{2} \) on all intervals \( I_{l,J}\).}

\noindent  The proof is by induction on \( k \).

\noindent For \( k = 1 \) then by using lemma \ref{main_lemma} we
have that for \( J \) and \( l \) big enough
\[
 \frac{1}{p_1(u_{l,J})} \sum_{n=1}^{p_1(u_{l,J})} \frac{1}{J}
\sum_{j=1}^J b_{u_{l,J}+j} 1_A(n) 1_A(p_1(u_{l,J}+j) - n) =
\]
\[
<1_A,\frac{1}{J} \sum_{j=1}^J b_{u_{l,J}+j}(\xi(p_1(u_{l,J}+j) -
n) + d(A))>_{p_1(u_{l,J})} \,\, \geq
\]
\[ - \varepsilon +
d(A)<1_A,\frac{1}{J} \sum_{j=1}^J b_{u_{l,J}+j}>_{p_1(u_{l,J})} \,\,
> (1-\delta) d(A)^2 \frac{d^{*}(F)}{3}.\]

\noindent
Assume the statement of the claim holds for \( k \).  Let
\((p_1,\ldots,p_{k},p_{k+1})\) be polynomials of the same degree such that \(p_1 \) is the
 "biggest" among them (see conditions of lemma \ref{main_lemma}).
Without loss of generality we can assume that \( \min_{2 \leq i
\leq k+1} {c_i} = c_{k+1}\).  Then
for sufficiently large \( J\) and \(l\)
\[
\frac{1}{p_1(u_{l,J})} \sum_{n=1}^{p_1(u_{l,J})} \frac{1}{J}
\sum_{j=1}^J b_{u_{l,J}+j} 1_A(n) 1_A(p_1(u_{l,J}+j) - n) \ldots
\]
\[
1_A(p_k(u_{l,J}+j) - n) 1_A(p_{k+1}(u_{l,J}+j) - n) =
\]
\[
<1_A,\frac{1}{J} \sum_{j=1}^J b_{u_{l,J}+j} 1_A(p_1(u_{l,J}+j) -
n) \ldots
\]
\[ 1_A(p_k(u_{l,J}+j) - n) (\xi(p_{k+1}(u_{l,J}+j) - n) +
d(A)) >_{p_1(u_{l,J})}  -
\]
\[
d(A) \frac{1}{p_1(u_{l,J})} \sum_{n= p_{k+1}(u_{l,J})}^{p_1(u_{l,J})}
1_A(n) \frac{1}{J} \sum_{j=1}^J b_{u_{l,J}+j} 1_A(p_1(u_{l,J}+j) -  n) \ldots
  1_A(p_k(u_{l,J}+j) - n) =
\]
\[
d(A)<1_A,\frac{1}{J} \sum_{j=1}^J b_{u_{l,J}+j} 1_A(p_1(u_{l,J}+j)
- n) \ldots 1_A(p_k(u_{l,J}+j) - n)>_{p_1(u_{l,J})}+
\]
\[
<1_A,\frac{1}{J} \sum_{j=1}^J b_{u_{l,J}+j}1_A(p_1(u_{l,J}+j) - n)
\ldots 1_A(p_k(u_{l,J}+j) - n) \xi(p_{k+1}(u_{l,J}+j) -
n)>_{p_1(u_{l,J})}  -
\]
\[
 d(A) \frac{1}{p_1(u_{l,J})} \sum_{n= p_{k+1}(u_{l,J})}^{p_1(u_{l,J})}
1_A(n) \frac{1}{J} \sum_{j=1}^J b_{u_{l,J}+j} 1_A(p_1(u_{l,J}+j) -  n) \ldots
  1_A(p_k(u_{l,J}+j) - n)
\]
\[
  >
\]
\[
   d(A) \frac{1}{p_1(u_{l,J})} \sum_{n= 1}^{p_{k+1}(u_{l,J})}
1_A(n) \frac{1}{J} \sum_{j=1}^J b_{u_{l,J}+j} 1_A(p_1(u_{l,J}+j) -  n) \ldots
  1_A(p_k(u_{l,J}+j) - n)
   - \varepsilon
\]
\[
  > d(A)c(1-\delta')d(A)^{k+1}\frac{d^*(F)}{3}
\]
\[
  > c(1 - \delta)d(A)^{k+2}\frac{d^{*}(F)}{3}.
\]
We used claim \( 2 \) in the first inequality and induction hypothesis in
the second inequality. This ends the proof of claim \( 3 \).

%

 \noindent  By the
definition of \( F \)  it follows that for every non-zero value of
\[
a_{u_{l,J}+j} 1_A(n) 1_A(p_1(u_{l,J}+j) - n) 1_A(p_2(u_{l,J}+j) -
n) \ldots 1_A(p_{k-1}(u_{l,J}+j)-n) \]
 (thus it equals to one), the remaining factor in the summands of \( B_{u_{l,J},J} \)
 is negative, namely, \( \xi(p_k(u_{l,J}+j)-n) = -d(A)
\). Therefore, by using claim \( 3 \) we get \( |B_{u_{l,J},J}| \geq c (1 - \varepsilon) d^{k+1}(A)
\frac{d^{*}(F)}{3} \) for any \( l \) and for \( J \) big
enough. Thus \( |B_{u_{l,J},J}| \) is bounded from zero.
\newline
On the other hand, by claim \( 2 \)  it follows
that for any \( \varepsilon > 0 \) there exists \( J= J(\varepsilon)
\) and \( N = N(J(\varepsilon))\) such that \( |B_{N,J}| <
\varepsilon \). Therefore we get a contradiction.

\noindent    We have proved that the set of all \( z\)'s such that the
statement of the theorem holds has a lower density one. Therefore it
intersects every set of positive density (even of positive upper
density), in particular, \(A\).

\hspace{12cm} \qed
\end{proof}

\newpage
\section{The equation
$xy=z$ and normal sets}\label{normal_section} 
\subsection{Normal sets and diophantine
equations}\label{norm_sets_Diophantine _eq}
 We expect that there are
many diophantine equations which are solvable in every normal set.
 We denote by \( DSN \) the
family of diophantine equations (including systems of
equations) which are solvable within every normal set.
\newline
It is easily seen that the equation \( x+y = z \) is in \( DSN \).
This equation is called the additive Schur equation. Schur proved that
the equation is "partition regular". This means that for any finite coloring
 of \( \mathbb{N} \), there exists a monochromatic solution for the
 equation (see
\cite{schur}).
\newline
 Systems of linear diophantine equations that are partition regular
  are classified by Rado in \cite{rado}.
 Such systems are usually called Rado systems.
\newline
In section \ref{sub_sect_Rado} we show that any Rado system of linear equations is in \( DSN \).

\noindent At the moment, we don't know the richness of the \( DSN \) family. By
the aforementioned result, a large family of linear equations
(Rado's systems) are in \( DSN \). For non-linear case, we don't
know much. For instance, it is not known whether the equation \( x^2
+ y^2 = z^2 \) is in \( DSN \). In this chapter we prove
that the equation \( xy = z \) is not in \( DSN \). The last
equation is called the multiplicative Schur equation. It should be
mentioned that for partitions of \( \mathbb{N} \) into finite number
of subsets, at least one of subsets contains solutions for both
Schur's additive and multiplicative equations (see \cite{berg1}).
Therefore, there exist partition regular equations that
are not in \( DSN \). We use the notion of Liouville's function to
construct a normal set in which the multiplicative Schur's equation
is unsolvable.
\newline
\begin{definition}
Liouville's function \( \lambda : \mathbb{N} \rightarrow \{ -1, 1 \}
\) is defined as follows:
\[
\lambda(p_1^{e_1}p_2^{e_2} \ldots p_k^{e_k}) = (-1)^{e_1+e_2+ \ldots
+e_k}
\]
where \( p_1, \ldots , p_k \) are primes.
\end{definition}
It is a well known and very deep question whether the set \( A = \{
n \in \mathbb{N} | \lambda(n) = -1 \} \) forms a normal set, see
\cite{liouv_1} and \cite{liouv_2}.  It seems that at present we are
far away from resolving this outstanding problem. But just for
clarity, if the answer for the question
is positive, then the aforementioned set \( A \) gives us an example
of a normal set with no solution to the equation \( xy=z \).
\newline
For the following we will use a modified  Liouville's function \(
\lambda_Q \) which is defined by random choice of subset \( Q \)
inside \(P\) (prime numbers) as follows
\[
\lambda_Q(p_1^{e_1}p_2^{e_2}\ldots p_k^{e_k}) = \lambda_Q(p_1)^{e_1}
\lambda_Q(p_2)^{e_2} \ldots \lambda_Q(p_k)^{e_k}
\]
and
\begin{displaymath}
\lambda_Q(p) = \left\{
                       \begin{array}{cc} -1 &  p \in Q \\
                                   1 &  p \not \in Q
                                   \end{array}
                                   \right.
\end{displaymath}
 By randomness of \( Q \) we mean  that a choice of
every prime number \( p \)  is independent of other prime numbers
and \( Pr ( p \in Q ) = 0.5 \)  for any \( p \in P \).
\newline
One defines \(A_Q = \{ n \in \mathbb{N} | \lambda_Q(n) = -1 \}\).  We prove  the following
\begin{theorem}
\label{main_theorem} For almost every \( Q \) the set \( A_Q \) is
normal.
\end{theorem}
This theorem gives us an infinite family of normal sets such that
the multiplicative Schur's equation is not solvable in these sets.
\newline
In the section (\ref{xy=z^2}) we prove that the equations \( xy=z^2
\), \( x^2 + y^2 = square\) and \( u^2 - v^2 = square \) are in
\(DSN \).

\subsection{ $A_Q$ is normal for a.e. $Q$}\label{{A_Q_is_normal}}
 We start from an obvious claim about
normality of \( A_Q \) which is a restatement of lemma \ref{norm_lemma}.
\begin{lemma}
\label{norm_lemma_2} Let \( Q \subset P \) be given, then \( A_Q \) is
a normal set \( \Leftrightarrow \) for any \(k \in (\mathbb{N} \cup
\{0\}) \) and any \( i_1 < i_2 <  \ldots < i_k \) we have
\[
 \lim_{N \rightarrow \infty} \frac{1}{N} \sum_{n=1}^N \lambda_Q(n) \lambda_Q(n+i_1) \ldots \lambda_Q(n+i_k)
 =0.
\]
\end{lemma}


\noindent
Denote
\begin{equation}
\label{T_N}
  T_N = \sum_{n=1}^{N} \lambda_Q(n) \lambda_Q(n+i_1) \ldots \lambda_Q(n+i_k).
\end{equation}
The next step is to show
\begin{displaymath}
 \sum_{N=1}^{\infty} E((\frac{1}{N^{40}} \sum_{n=1}^{N^{40}} \lambda_Q(n) \lambda_Q(n+i_1) \ldots \lambda_Q(n+i_k))^2) <
 \infty.
\end{displaymath}
\begin{lemma}
\label{main_l}
 With \( T_N  \) as  defined in (\ref{T_N}),
\( E( T_N^2 ) \leq O(\frac{1}{N^{0.05}}) \).
\end{lemma}
\begin{proof}
By linearity of expectation we get
\[
 E(T_N^2) = \frac{1}{N^2} \sum_{x,y=1}^{N} E(\lambda_Q(x) \lambda_Q(x+i_1) \ldots \lambda_Q(x+i_k)
\lambda_Q(y) \lambda_Q(y+i_1) \ldots \lambda_Q(y+i_k)).
\]
Note that for any \( m \in \mathbb{N} \), \( E(\lambda_Q(m)) = 0 \)
unless \( m \) is a square in which case \( E(\lambda_Q(m)) = 1 \).
\newline
Let us denote by
\[ \phi(x) \risingdotseq \lambda_Q(x) \lambda_Q(x+i_1) \ldots
\lambda_Q(x+i_k) \]
 and
 \[ \xi(x) \risingdotseq x(x+i_1) \ldots
(x+i_k). \]
 By distribution of \( Q \) we get
 \[ E (\phi(x)
\phi(y)) = 1    \Leftrightarrow  \xi(x) \xi(y) = m^2. \]
 Otherwise
\[ E (\phi(x) \phi(y)) = 0. \]
 Therefore, to obtain  an upper bound
on \( E({T_N}^2) \), we give an upper bound on the number of pairs
\( (x,y) \in [1,N] \times [1,N] \) which satisfy \( \xi (x) \xi(y) =
square \).
\newline
For a given \( x \in [1,N] \) let us assume that \( \xi(x) = c_x m^2
\), where \( c_x \) is a square-free number, say \( c_x = p_{j_1}
\ldots p_{j_l} \) is the prime factorization of \( c_x \). Then we
will define \( h(x) = l \) (thus \( h(x) \) is a number of primes in
prime factorization of maximal square-free number which divides \( x
\)).
Denote by \( D \) the set of all possible common divisors  of the
numbers \( x, x+i_1, \ldots, x+i_k \) (i.e. positive integers which
divide at least two of them).  For a finite non empty set \( S \) of
positive numbers we denote by \( m(S) \) the product of all elements
of \( S \)  and, for empty set, we fix  \( m(\emptyset) = 1 \).
\newline
 Note that \( \xi(x) \xi(y) = square \) \( \Rightarrow \)
there exist \( S_1 \subset D \) and \( S_2 \subset \{p_{j_1}, \ldots
,p_{j_l} \}\) such that \( y = m(S_1) m(S_2) square \).
\newline
Assume  \( |D| = r \) (\( r \) depends only on the set \( \{i_1,
\ldots, i_k \} \) and doesn't depend on \(x\)).
Then we obtain \( \xi(x) \xi(y) = square \) for at most \(
2^r2^{h(x)} \sqrt{N} \) \( y \)'s inside \( [1,N] \).
 Thus
\[
E(T_N^2) \leq \frac{1}{N^2} ( \sum_{n=1}^N 2^r2^{h(n)} \sqrt{N} )
\leq  \frac{c}{N^{1.5}} \sum_{n=1}^N 2^{h(n)}
\]
Therefore it remains to bound the expression \( \sum_{n=1}^N
2^{h(n)} \).
\newline
If \( \xi(n) \) does not contain as dividers \( 2,3 \) then \( h(n)
\leq \log_5{(n+i_k)^{k+1}} = (k+1) \frac{\log_2{(n+i_k)}}{\log_2{5}}
\). This gives us
\[
2^{h(n)} \leq 2^{k+1} (n+i_k)^{\frac{1}{\log_2{5}}} \leq C_1
(n+i_k)^{0.45}
\]
But if \( \xi(n) \) contains \( 2\) or \( 3\) as dividers then \(
h(n) \) can increase by at most two, this means \( 2^{h(n)} \leq 2^2
C_1 (n+i_k)^{0.45} \). Thus \( \sum_{n=1}^N 2^{h(n)} \leq C_2
(N+i_k)^{1.45} \) and therefore we get
\[
E(T_N^2) \leq C_3 \frac{1}{N^{0.05}}.
\]

\hspace{12cm} \qed
\end{proof}

\noindent
\begin{proof}(\textbf{theorem \ref{main_theorem})}
From the last lemma we conclude that \(  \sum_{N=1}^{\infty}
E(T_{N^{40}}^2) < \infty \). Thus almost surely
\( T_{N^{40}} \rightarrow 0\). By lemma \ref{tec_lemma}  it follows that almost
surely \( T_N \rightarrow 0\). And from lemma \ref{norm_lemma_2}
(and countability of necessary conditions) it follows that for
almost all \( Q \subset P \) the sets \( A_Q \) are normal.

\hspace{12cm} \qed
\end{proof}

\noindent
We can now demonstrate the main result of this note.
\begin{theorem}
There exists  a normal set \( A \subset \mathbb{N} \) such that the
multiplicative Schur's equation is not solvable inside \( A \).
\end{theorem}
\begin{proof}
 We have already shown the existence of many \(Q \) ( \( Q \subset P \))  such that \( A_Q \) are normal.
By definition of \( A_Q \) follows that for any \( x,y \in A_Q \)
the number \( xy \not \in A_Q \). Therefore we can't find \( x,y,z
\in A_Q \) such that \( xy = z \).

\hspace{12cm} \qed
\end{proof}
\begin{corollary}
For any equation \( xy = c n^k \) (where \( c,k \) are natural
numbers, \( c \) is not a square and \( k \) is even) we can find a
normal set \( A_{c,k} \subset \mathbb{N} \) such that for any \( x,y
\in A \) we have \( xy \not = c n^k \) for every natural \( n \).
\end{corollary}
\begin{proof}
We take \( A_Q \)  be a normal and such that \( \lambda_Q(c) = -1 \)
(it happens with the positive probability \(\frac{1}{2} \), and thus
there exist such sets). Then obviously we can't solve the proposed
equation inside \( A_Q \).

\hspace{12cm} \qed
\end{proof}

\subsection{ Solvability of the equation $xy = z^2$ and related
problems}\label{xy=z^2}
\begin{theorem}
\label{xy=z^2}
(theorem \ref{th1_pos} of \(\S\)\ref{normal_sets_xyz})
\newline
Let \( A \subset \mathbb{N} \) be a WM set. Then
there exist \( x,y,z \in A \) (\( x \not = y \)) such that \( xy =
z^2 \).
\end{theorem}
\begin{proof}
For a set \( S \subset \mathbb{N} \) let us define \( S_a = \{ n
\in \mathbb{N} | an \in S \} \), where \( a \in \mathbb{N} \). It
is easily seen that if \( S \) is a WM set then \( S_a \) is again
a WM set with the same statistics as \( S \) for any natural \( a
\) (see \cite{furst2}). We denote by \( d(S) \) density of a set
\( S \), if it exists.
\newline
Let \( A \) be a WM set. We denote by \( R_n \risingdotseq A_{2^n}
\). For any \( n \), \( d(R_n) = \frac{1}{2} \). Let us
denote by
\[\mu_N (S) = \frac{|S \cap \{ 1, 2, , \ldots N
\}|}{N} \]
 for any \( S \subset \mathbb{N} \) and any \( N \in
\mathbb{N} \).
\newline
By Szemer\'{e}di's theorem (finite version), for any \( \delta > 0\)
and any \( l \in \mathbb{N} \) there exists \( N(l,\delta) \) such
that for any \( N \geq N(l,\delta) \) and any \( F \subset
\{1,2,\ldots,N\} \) such that \( \frac{|F|}{N} \geq \delta \) the
set \( F \) contains an arithmetic progression of length \( l \)
(see \cite{szemeredi}).
\newline
One chooses \( K \geq  N(3,\frac{1}{3}) \). Then there exists \( N_K
\) such that \( \mu_{N_K} (R_i) \geq \frac{1}{3} \) for every \( 1
\leq i \leq K \).
\newline
We claim that there exists \( F \subset \{1,2, \ldots, K \} \) such
that \( \frac{|F|}{K} \geq \frac{1}{3} \) and \( \mu_{N_K}(\cap_{j
\in F} R_j) > 0 \). If not, let us denote \( 1_{R_i} \) to be the
indicator function of the set \( R_i \) inside the set \( \{ 1,
\ldots , N_K \} \). Then \[ \int_{[1,N_K]} (1_{R_1} + \ldots +
1_{R_K}) d\mu_{N_K}  = \sum_{j=1}^{K} \int_{[1,N_K]} 1_{R_j}
d\mu_{N_K}  \geq \frac{K}{3}. \]
\newline
Therefore  \( \exists n: \, 1 \leq n \leq N_K \) such that
\( \sum_{j=1}^K 1_{R_j}(n) \geq \frac{K}{3} \).
\newline
Thus  \( \mu_{N_K}(\cap_{j
\in F} R_j) > 0 \).
\newline
Let \( F \subset \{1,2, \ldots, K \} \) such that \( \frac{|F|}{K}
\geq \frac{1}{3} \) and \( \mu_{N_K}(\cap_{j \in F} R_j) > 0 \).
Then by the choice of \( K \) it follows that \( F \) necessarily
contains arithmetic progression of length \( 3 \). The last statement means
there exist \( a,b,c \in F \) such that \( a+c=2b \). Let us take \(
R_a, R_b , R_c \). We have \( R_a \cap R_b \cap R_c \neq \emptyset
\) and this means there exists \( n \in \mathbb{N} \) such that \(
n2^a \in A \) and \( n2^b \in A \) and \( n2^c \in A \). Let us
denote by \( x,y,z \) the following elements of \( A \):
 \( x = n2^a \), \( y = n2^c \), \( z = n2^b \). Then we have
\[
xy = z^2.
\]

\hspace{12cm} \qed
\end{proof}

\noindent
\textbf{Question:} Are the equations \( x y = c^2 z^2 \), where \( c
> 0 \) is a natural number, always solvable inside an arbitrary
normal set?

\noindent We repeat the formulation of theorem
\ref{th2_pos}.
\nonumber \begin{theorem}
 Let \( A \subset \mathbb{N} \) be an arbitrary
normal set. Then there exist \( x,y,u,v \in A \) such that \( x^2 +
y^2 = square \) and \( u^2 - v^2 = square \).
\end{theorem}
\begin{proof}
Note that there exist \( a,b,c \in \mathbb{N} \) such that \( a^2 +
b^2 = square \) and \( a^2 + c^2 = square \) and \( b^2 + c^2 =
square \). For example \( a = 44, b = 117, c = 240 \).
\newline
Let \( A \subset \mathbb{N} \) be an arbitrary normal set. We look
at \( A_a, A_b, A_c \) which are defined as in the proof of theorem \ref{xy=z^2}. Then
\( d(A_a) = d(A_b) = d(A_c) = \frac{1}{2} \) and thus it can not be
true that the intersection of each pair from the triple is empty.
\newline
Without loss of generality, let us assume that \( A_a \cap A_b \neq
\emptyset \).
\newline
Thus there exists \( z \in A_a \cap A_b \) or equivalently \( z a, z
b  \in A \). But \( a^2 + b^2 = square \) and therefore \( (za)^2 +
(zb)^2 = square \).
\newline
A proof that the equation \(  u^2 - v^2 = square \) is solvable in
any normal set is similar. We use the fact that there exist \( a, b
,c \in \mathbb{N} \) with \( a < b < c \) such that \( c^2 -
b^2 = square \) and \( c^2 - a^2 = square \) and
 \(b^2 - a^2 = square \). For example \( a = 153, b = 185, c = 697 \).

\hspace{12cm} \qed
\end{proof}

\noindent \textbf{Question:} For an arbitrary normal set \( A \) do
there exist \( x,y,z \in A \) such that \( x^2 + y^2 = z^2\)?
\newpage
\section{Appendix}\label{appendix}
\numberwithin{lemma}{section} \numberwithin{theorem}{section}
\numberwithin{prop}{section} \numberwithin{remark}{section}
\numberwithin{definition}{section} \numberwithin{corollary}{section}

\noindent In this section we prove all technical lemmas and
propositions that were used in the thesis.
\newline
We start with the key lemma which is a finite modification of
Bergelson's lemma in \cite{berg_pet} and its origin is in lemma of
van der Corput.
\begin{lemma}\( \rm( \)\textit{van der Corput}\( \rm) \)
\label{vdrCorput}
 Suppose \(\varepsilon
>0 \) and \( \{u_{j}\}_{j=1}^{\infty } \) is a family
of vectors in  Hilbert space, such that \( \Vert u_j \Vert \leq 1 \,
\rm( 1 \leq j \leq \infty \rm)  \). Then there exists \(
I'(\varepsilon) \in \mathbb{N} \), such that for every \( I \geq
I'(\varepsilon) \) there exists \( J'(I,\varepsilon) \in \mathbb{N}
\), such that the following holds:
\newline
For \( J \geq J'(I,\varepsilon) \) for which we obtain
\[
\left| \frac{1}{J}\sum ^{J}_{j=1}<u_{j}, u_{j+i}>\right| <
\frac{\varepsilon}{2}, \]  for set of \( i \)'s in the interval \(
\{1,\ldots,I\} \) of density \( 1 - \frac{\varepsilon}{3}\) we have
\[
\left\Vert \frac{1}{J}\sum _{j=1}^{J}u_{j}\right\Vert < \varepsilon.
\]
\end{lemma}
\begin{proof}
For an arbitrary \( J \) define \( u_{k}=0 \) for every \(
\textrm{k}<1 \) or \( k>J \). The following is an elementary
identity:
\[
\sum ^{I}_{i=1}\sum ^{J+I}_{j=1}u_{j-i}=I\sum
^{J}_{j=1}u_{j}.
\]
Therefore,  the inequality \( \left\Vert \sum
_{i=1}^{N}u_{i}\right\Vert ^{2}\leq N\sum _{i=1}^{N}\left\Vert
u_{i}\right\Vert ^{2} \) yields
 \[
\left\Vert I\sum ^{J}_{j=1}u_{j}\right\Vert ^{2}\leq (J+I)\sum
^{J+I}_{j=1}\left\Vert \sum ^{I}_{i=1}u_{j-i}\right\Vert ^{2}=\]
 \[
(J+I)\sum ^{J+I}_{j=1}<\sum ^{I}_{p=1}u_{j-p},\sum
^{I}_{s=1}u_{j-s}>=\]
 \[
(J+I)\sum ^{J+I}_{j=1}\sum ^{I}_{p=1}\left\Vert u_{j-p}\right\Vert
^{2}+2(J+I)\sum ^{J+I}_{j=1}\sum
^{I}_{r,s=1;s<r}<u_{j-r},u_{j-s}>=\]
 \[
 (J+I)(\Sigma _{1}+2\Sigma _{2}),
\]
where \( \Sigma _{1}=I\sum ^{J}_{j=1}\left\Vert u_{j}\right\Vert
^{2} \) by  the aforementioned elementary identity and
 \(
\Sigma _{2}=\sum ^{I-1}_{h=1}(I-h)\sum ^{J}_{j=1}<u_{j},u_{j+h}> \).
The last expression is obtained by rewriting \( \Sigma_{2} \), where
\( h=r-s \). By dividing the foregoing inequality  by \( I^{2}J^{2}
\) we obtain
 \[
\left\Vert \frac{1}{J}\sum ^{J}_{j=1}u_{j}\right\Vert
^{2}<\frac{J+I}{IJ}+\frac{J+I}{J}\left( \frac{\varepsilon}{2} +
\frac{\varepsilon}{3} \right) = \frac{J+I}{J}\left( \frac{1}{I} +
\frac{5 \varepsilon}{6} \right).
\]
Choose \( I'(\varepsilon) \in \mathbb{N} \), such that
\(\frac{12}{\varepsilon}
 \leq I'(\varepsilon) \leq \frac{12}{\varepsilon}+1\). Then for
every \( I \geq I'(\varepsilon) \) we have \(  \frac{1}{I} + \frac{5
\varepsilon}{6} \leq \frac{11\varepsilon}{12}\). There exists \(
J'(I,\varepsilon) \in \mathbb{N} \), such that for every \( J \geq
J'(I,\varepsilon) \) we obtain \( \frac{J+I}{J} < \frac{12}{11} \).
As a result, for every \( I \geq I'(\varepsilon) \) there exists \(
 J'(I,\varepsilon) \), such that for every \( J \geq
 J'(I,\varepsilon) \)
\[
\left\Vert \frac{1}{J}\sum ^{J}_{j=1}u_{j}\right\Vert ^{2} <
\varepsilon.
\]

\hspace{12cm} \qed
\end{proof}

\noindent The next proposition is useful in section
\ref{lin_proof_suff}.
\begin{prop}
\label{wm_prop}
 Let \( A \subset \mathbb{N} \) be a WM-set. Then
for every integer \( a>0 \) and every integers \( b_1,b_2,\ldots,b_k
\) we obtain the following
\begin{displaymath}
 \lim_{N \rightarrow \infty} \frac{1}{N}
\sum_{n=1}^N \xi(n+b_1)\xi(n+b_2) \ldots \xi(n+b_k) =
\end{displaymath}
\[
 \lim_{N \rightarrow \infty} \frac{1}{N} \sum_{n=1}^N
\xi(an+b_1)\xi(an+b_2) \ldots \xi(an+b_k),
\]
where \( \xi \doteq 1_A - \dd(A) \).
\end{prop}
\begin{proof}
Consider the weak-mixing measure preserving system \( ( X_{\xi},
\mathbb{B},\mu,T) \).
\newline
The left side of the equation in the proposition is \( \int_{X_{\xi}}
T^{b_1} f T^{b_2} f \ldots T^{b_k} f d\mu \), where \( f(\omega)
\doteq \omega_0 \) for every infinite sequence inside \( X_{\xi}\).
We make use of the notion of disjointness of measure preserving systems. By \cite{furst2} we know that
every weak-mixing system is disjoint from any Kronecker system which
is a compact monothethic  group with Borel \(\sigma\)-algebra, the
Haar probability measure, and the shift by an a priori chosen
element of the group. In particular, every weak-mixing system is
disjoint from the measure preserving system \(
(\mathbb{Z}_a,\mathbb{B}_{\mathbb{Z}_a}, S , \nu)\), where
 \( \mathbb{Z}_a = \mathbb{Z}/ a\mathbb{Z} \), \( S(n)
\doteq n+1 (\mod a) \). The measure and the \( \sigma\)-algebra of the last
system are uniquely determined. Therefore, from  Furstenberg's
theorem
(see \cite{furst2}) it follows that the point \( (\xi,0) \in X_{\xi}
\times \mathbb{Z}_a \) is a generic point of the product system \(
(X_{\xi} \times \mathbb{Z}_a, \mathbb{B} \times
\mathbb{B}_{\mathbb{Z}_a}, T \times S, \mu \times \nu) \). Thus, for
every continuous function \( g \) on \( X_{\xi} \times \mathbb{Z}_a
\) we obtain
\[
\int_{X_{\xi} \times \mathbb{Z}_a} g(x,m) d\mu(x) d \nu(m) = \lim_{N
\rightarrow \infty}  \frac{1}{N} \sum_{n=1}^N g(T^n \xi, S^n 0).
\]
Let \( g (x,m) \doteq f(x)1_0(m) \) which is obviously continuous on
\( X_{\xi} \times \mathbb{Z}_a \). Then  genericity of the point \(
(\xi,0) \) yields
\[
\int_{X_{\xi} \times \mathbb{Z}_a} f(x)1_0(m) d\mu(x) d \nu(m) =
\frac{1}{a} \int_{X_{\xi}} f(x) d \mu(x) = \]
\[
\lim_{N \rightarrow \infty}  \frac{1}{N} \sum_{n=1}^N f(T^n \xi)
1_0(n) = \lim_{N \rightarrow \infty}  \frac{1}{a} \frac{1}{N}
\sum_{n=1}^N f(T^{an} \xi).
\]
Taking instead of the function \( f \) the continuous function \(
T^{b_1}f T^{b_2} f \ldots T^{b_k} f\)  in the definition of \( g \)
finishes the proof.

\hspace{12cm} \qed
\end{proof}

\noindent The following lemma is simple fact that for a weak-mixing
system \( X \) not only an average of shifts for a function converge
to a constant in \( L^2 \) norm but also  weighted averages (weights
are bounded) converge to the same constant.
\begin{lemma}
\label{wm_sub_lemma1} Let \( ( X, \mathbb{B}, \mu, T ) \) be a
weak-mixing system  and \( f \in L^{2}(X) \) with \( \int_X f d\mu =
0 \). Let \( \varepsilon > 0 \). Then there exists \( \mathbb{J}
> 0 \) such that for any \( J
> \mathbb{J} \) we have
\[
 \left\| \frac{1}{J} \sum_{j=1}^J b_j T^{j} f \right\|_{L^2(X)} <
 \varepsilon
\]
for any sequence \(  b = (b_1,b_2, \ldots , b_n, \ldots) \in \{ 0 ,1
\}^\mathbb{N}  \).
\end{lemma}
\begin{proof}
Let \( \varepsilon > 0 \).
\newline
 By one of the properties of weak
mixing, for any \( f \in L^2(X) \) with \( \int_X f d \mu(x) = 0 \)
we have \( \frac{1}{N} \sum_{n=1}^N |<T^n f, f>| \rightarrow 0 \).
\newline
We denote by \( c_n = c_{(-n)}= |<T^n f, f>| \) and we have that \(
\frac{1}{N} \sum_{n=1}^N c_n \rightarrow 0 \). Then for any \(
\varepsilon > 0 \) there exists \( \mathbb{J} > 0 \) such that for
any \( J
> \mathbb{J} \) we have
\[
\left\| \frac{1}{J} \sum_{j=1}^J b_j T^{j} f \right\|^2 \leq
\frac{1}{J^2} \sum_{j=1,k=1}^J b_j b_k c_{j-k} \leq \frac{1}{J^2}
\sum_{j=1,k=1}^J c_{j-k} \leq \varepsilon.
\]

\hspace{12cm} \qed
\end{proof}

\noindent The next two lemmas are very useful for constructing
normal sets with specifical properties (we use them in this thesis for
constructing counterexamples).
\begin{lemma}
\label{norm_lemma} Let \( A \subset \mathbb{N} \). Let \( \lambda(n)
= 2 1_A(n) - 1 \). Then \( A \) is a normal set \( \Leftrightarrow \)
for any \(k \in (\mathbb{N} \cup \{0\}) \) and any \( i_1 < i_2 <
\ldots < i_k \) we have
\[
 \lim_{N \rightarrow \infty} \frac{1}{N} \sum_{n=1}^N \lambda(n) \lambda(n+i_1) \ldots \lambda(n+i_k)
 =0.
\]
\end{lemma}
\begin{proof}
"\(\Rightarrow\)" If \( A \) is normal then any finite word \( w \in
\{ -1, 1\}^{*} \)
 has the "right" frequency \( \frac{1}{2^{|w|}}\) inside
\( w_{A} \). This guarantees that "half of the time" the function \(
\lambda(n) \lambda(n+i_1) \ldots \lambda(n+i_k) \) equals \(1 \) and
"half of the time" is equal to \( -1\).
 Therefore we get the desired conclusion.
\newline
"\( \Leftarrow \)"  Let \( w \) be an arbitrary finite word of plus
and minus ones: \( w = a_1 a_2 \ldots a_k \) and we have to prove
that \( w \) occurs in \( w_{A}\) with the frequency \( 2^{-k} \).
For every \( n \in \mathbb{N} \) the word \( w \) occurs in \( 1_A \) and starting from \( n \)
if and only if
\[
 \left\{ \begin{array}{lll} 1_A(n) = a_1 \\
 \ldots \\
1_A(n+k-1) = a_k
\end{array}
\right.
\]
The latter is equivalent to the following
\[
 \left\{ \begin{array}{lll} \lambda(n) = 2a_1 - 1 \\
 \ldots \\
\lambda(n+k-1) = 2 a_k -1
\end{array}
\right.
\]
The frequency of \( w \) within \( 1_A \) is equal to
\[
  \lim_{N \rightarrow \infty} \frac{1}{N} \sum_{n=1}^N \frac{\lambda(n) (2 a_1 -1) +1}{2} \ldots
  \frac{\lambda(n+k-1) (2 a_k - 1) + 1}{2}.
\]
By assumptions of the lemma the latter expression is equal to \( \frac{1}{2^k}\).

\hspace{12cm} \qed
\end{proof}

\begin{lemma}
\label{tec_lemma} Let \( \{ a_n \} \) be a bounded sequence. Denote
by \( T_N = \frac{1}{N} \sum_{n=1}^N a_n \). Then \( T_N \)
converges to a limit \( t \) \( \Leftrightarrow \) there exists a
sequence of increasing indices \( \{ N_i \} \) such that \(
\frac{N_i}{N_{i+1}} \rightarrow 1 \) and \( T_{N_i} \rightarrow_{i
\rightarrow \infty} t \).
\end{lemma}

\newpage

\newpage

\end{document}